\documentclass[12pt,a4paper]{amsart}
\usepackage{amsfonts}
\usepackage{amsthm}
\usepackage{amsmath}
\usepackage{amscd}
\usepackage[latin2]{inputenc}
\usepackage{t1enc}
\usepackage[mathscr]{eucal}
\usepackage{indentfirst}
\usepackage{graphicx}
\usepackage{graphics}
\usepackage{pict2e}
\usepackage{epic}
\numberwithin{equation}{section}
\usepackage[margin=2.9cm]{geometry}
\usepackage{epstopdf} 
\usepackage{enumerate}
\usepackage{verbatim}
\usepackage{cite}
\usepackage{mathtools}
\usepackage{mathrsfs}
\usepackage{xcolor}
\usepackage{hyperref}
\usepackage{relsize}
\usepackage{esint}
\usepackage{comment}

\theoremstyle{plain}
\newtheorem{Th}{Theorem}[section]
\newtheorem{Lemma}[Th]{Lemma}
\newtheorem{Cor}[Th]{Corollary}
\newtheorem{Prop}[Th]{Proposition}

\newtheoremstyle{named}{}{}{\itshape}{}{\bfseries}{.}{.5em}{\thmnote{#3}}
\theoremstyle{named}
\newtheorem*{namedtheorem}{Theorem}

\theoremstyle{definition}

\newtheorem{Conj}[Th]{Conjecture}

\definecolor{ao(english)}{rgb}{0.0, 0.5, 0.0}

\begin{document}

\title[Average area ratio and normalized total scalar curvature]{Average area ratio and normalized total scalar curvature of hyperbolic n-manifolds}

\author{Ruojing Jiang}

\address{the University of Chicago, Department of Mathematics, Chicago, IL 60637} 
\email{ruojing@math.uchicago.edu}

 \subjclass[2019]{} 
 \date{}

\begin{abstract} 
On closed hyperbolic manifolds of dimension $n\geq 3$, we review the definition of the average area ratio of a metric $h$ with $R_h\geq -n(n-1)$ relative to the hyperbolic metric $h_0$, and we prove that it attains the local minimum of one at $h_0$, which solves a local version of Gromov's conjecture. Additionally, we discuss the relation between the average area ratio and normalized total scalar curvature for hyperbolic $n$-manifolds, as well as its relation to the minimal surface entropy if $n$ is odd.
\end{abstract}

\maketitle

\section{Introduction}
On a closed negatively curved manifold $M$, a useful fact is that the unit tangent bundle admits a one-dimensional foliation whose leaves are orbits of the geodesic flow. Given a pair of Riemannian metrics $h$ and $h_0$ on $M$, the comparisons of the geometric objects associated with $h$ and $h_0$
have been studied and become well-known. 
For instance, the geodesic stretch $I_\mu(h/h_0)$ measures the stretching of the metric $h$ relative to the reference metric $h_0$ and the measure $\mu$. Comparisons related to $I_\mu(h/h_0)$ were discussed by Knieper \cite{Knieper}.
Another example is the volume entropy $E_{vol}$ that counts the number of closed geodesics, the locally symmetric metric attains the minimum among all metrics on $M$ with the same volume (see Besson-Courtois-Gallot \cite{Besson-Courtois-Gallot}). The more detailed results will be reviewed later on, for the light they shed on what could be expected and extended to dimensions larger than one.

The first analogue to dimension two was discussed by Gromov in \cite{Gromov}. A two-dimensional foliation of $Gr_2M$ is formed by a family of stable minimal surfaces of $M$. In particular, let $h_0$ denote the hyperbolic metric, there is a canonical foliation of $Gr_2M$ whose leaves are totally geodesic planes. Gromov introduced the average area ratio $\text{Area}_F(h/h_0)$ which is similar in spirit to the geodesic stretch, to measure the ``stretching'' of the area of leaves on $Gr_2M$ of the metric $h$ relative to the hyperbolic metric $h_0$ under the map $F$. The comparison result and the interpretation of $\text{Area}_{\text{Id}}(h/h_0)$ for $3$-manifolds were studied by Lowe-Neves \cite{LN2021} using Ricci flow. In this paper, we focus on the local picture of $\text{Area}_{\text{Id}}(h/h_0)$ for any dimension $n\geq 3$, and investigate a more general case for $3$-manifolds.

On the other hand, Calegari-Marques-Neves \cite{CMN} defined the minimal surface entropy $E(h)$ based on the construction of surface subgroups by Kahn-Markovic \cite{Kahn-Markovic}, and the equidistribution property of $\text{PSL}(2,\mathbb{R})$-action on the space of minimal laminations (initiated by Ratner \cite{Ratner} and Shah \cite{Shah}, and recently formalized by Labourie \cite{Labourie}). $E(h)$ measures the number of essential minimal surfaces of $M$ with respect to $h$, thus shifting attention from a one-dimensional object (volume entropy) to dimension two. The minimal surface entropy of closed hyperbolic $3$-manifolds was computed by Calegari-Marques-Neves \cite{CMN} (also see Lowe-Neves \cite{LN2021}), and the author extended the result to manifolds of any odd dimension at least $3$ in \cite{Jiang}, we continue the discussion for higher dimensions in this paper.

\subsection{Gromov's conjecture}
Let $(M,h_0)$ be a closed hyperbolic $3$-manifold, and let $N$ be a closed $3$-manifold with Riemannian metric $g$ whose scalar curvature satisfies $R_g\geq -6$. $F: N\rightarrow  M$ is a smooth map with degree $d$. 
The average area ratio of $F$ is as follows. \begin{equation*}
    \text{Area}_F(g/h_0)=\int_{(y,P)\in Gr_2(M)}\sum_{x\in F^{-1}(y)}\lim_{\delta\rightarrow 0}\frac{\text{area}_g((dF_x)^{-1}(D_\delta))}{\delta} \,d\mu_{h_0},
\end{equation*}
where $D_\delta$ is a subset of the totally geodesic disc of $\mathbb{H}^n$ which is tangential to $P$ at $x$, $D_\delta$ has an area equal to $\delta$, and $\mu_{h_0}$ stands for the unit volume measure on $Gr_2(M)$ with respect to the metric induced by $h_0$. We refer to Section \ref{section_preliminaries} for a more detailed definition.
Gromov proved the following inequality using stability (see page 73 of \cite{Gromov}).  \begin{equation*}
    \text{Area}_F(g/h_0)\geq\frac{d}{3}.
\end{equation*}
The proof indicates that if the equality holds, then $(N,g)$ should be hyperbolic, and in the universal covers, the preimage of any totally geodesic disc is totally geodesic as well. However, it means the inequality is not sharp, and therefore, Gromov conjectured that the lower bound could be replaced by $d$. Moreover, the higher dimensional case may share the same property.

\begin{Conj}
Let $(M,h_0)$ be a closed hyperbolic manifold of dimension $n\geq 3$, and let $N$ be a closed $n$-manifold with Riemannian metric $g$, and its scalar curvature satisfies $R_g\geq -n(n-1)$. $F:N\rightarrow  M$ is a smooth map with degree $d$. We have
\begin{equation*}
    \text{Area}_F(g/h_0)\geq d.
\end{equation*}
The equality holds if and only if $F$ is a local isometry.
\end{Conj}

Recently, Lowe-Neves (Corollary 1.3 of \cite{LN2021}) verified the special case where $n=3$ and $F$ is a local diffeomorphism. We will generalize this conclusion in two aspects in Theorem \ref{Thm_gromov_high_codim} and Theorem \ref{Thm_gromov_3}.

\subsection{Minimal Surface Entropy}
Let $\mathbb{H}^n$ denote the hyperbolic $n$-space, where $n\geq 3$. In the Poincar\'{e} ball model, the asymptotic boundary $\partial_\infty \mathbb{H}^n$ can be considered as the $(n-1)$-unit sphere $S^{n-1}_\infty$. A homeomorphism $f:S^{n-1}_\infty\rightarrow S^{n-1}_\infty$ is called \emph{$K$-quasiconformal} if $K(f)=\underset{x\in S^{n-1}_\infty}{\text{ess sup}} K_f(x)\leq K$, where the \emph{dilatation} of $f$ is \begin{equation*}
        K_f(x)=\underset{r\rightarrow 0}{\lim\sup}\,\underset{y,z}{\sup}\frac{|f(x)-f(y)|}{|f(x)-f(z)|},
    \end{equation*}
    the supremum is taken over $y,z$ so that $r=|x-y|=|x-z|$.
A \emph{$K$-quasicircle} in $S^{n-1}_\infty$ is the image of a round circle under a $K$-quasiconformal map.

Let $M=\mathbb{H}^n/\pi_1(M)$ be a closed orientable $n$-manifold $(n\geq 3)$ that admits a hyperbolic metric $h_0$,
A closed surface immersed in $M$ with genus at least $2$ is said to be \emph{essential} if the immersion is $\pi_1$-injective, and the image of its fundamental group in $\pi_1(M)$ is called a \emph{surface subgroup}. 
Let $S(M,g)$ denote the set of surface subgroups of genus at most $g$ up to conjugacy, and let the subset $S(M,g,\epsilon)\subset S(M,g)$ consist of the conjugacy classes whose limit sets are $(1+\epsilon)$-quasicircles. Moreover, \begin{equation*}
    S_\epsilon(M)=\underset{g\geq 2}{\cup}S(M,g,\epsilon).
\end{equation*} 
Suppose $h$ is an arbitrary Riemannian metric on $M$. For any $\Pi\in S(M,g)$, we set \begin{equation*}
    \text{area}_h(\Pi)=\inf \{\text{area}_h(\Sigma):\Sigma\in\Pi\}.
\end{equation*}
Then the minimal surface entropy with respect to $h$ is defined as follows. 
\begin{equation}\label{eh}
    E(h)=\underset{\epsilon\rightarrow 0}{\lim}\,\underset{L\rightarrow\infty}{\lim\inf}\,\dfrac{\ln\#\{\text{area}_h(\Pi)\leq 4\pi(L-1):\Pi\in S_\epsilon(M)\}}{L\ln L}.
\end{equation}

According to \cite{CMN} and \cite{Jiang}, for odd-dimensional manifolds, and among metrics with sectional curvature less than or equal to $-1$, $E(h)$ attains the minimum at the hyperbolic metric $h_0$, and $E(h_0)=2$. On the other hand, Theorem 1.1 of \cite{LN2021} shows that when $n=3$, $E(h)$ is maximized at $h_0$ among all metrics with scalar curvature greater than or equal to $-6$.

\subsection{Main theorems}
In the following results, we assume $(M,h_0)$ is a closed hyperbolic manifold of dimension $n\geq 3$.
\begin{Th}\label{Thm_gromov_high_codim}
There exists a small neighborhood $\mathcal{U}$ of $h_0$ in the metric space of $M$, such that for any Riemannian metric $h\in\mathcal{U}$ on $M$ with $R_h\geq -n(n-1)$, we have 
 \begin{equation*}
    \text{Area}_{\text{Id}}(h/h_0)\geq 1.
    \end{equation*}
The equality holds if and only if $h=h_0$.
\end{Th}
Note that such a neighborhood is ``cylindrical'', that is, if $h\in\mathcal{U}$, then any metric $h'$ in the conformal class of $h$ with $R_{h'}\geq -n(n-1)$ also belongs to $\mathcal{U}$. 

Let $(N,g_0)$ be another closed hyperbolic manifold of the same dimension $n$, the theorem leads to the following corollary. 
\begin{Cor}
There exists a small neighborhood $\mathcal{U}$ of $g_0$ in the metric space of $N$, such that for any metric $g\in\mathcal{U}$ with $R_g\geq -n(n-1)$,
and for any local diffeomorphism $F: (N,g)\rightarrow (M,h_0)$ with degree $d$, we have \begin{equation*}
    \text{Area}_F(g/h_0)\geq d.
\end{equation*}
The equality holds if and only if $F$ is a local isometry between $g$ and $h_0$, i.e., $g=F^*(h_0)$.
\end{Cor}

\begin{proof}
Since $F$ is a local diffeomorphism, \begin{equation*}
    \text{Area}_F(g/h_0)=d\,\text{Area}_{\text{Id}}(g/F^*(h_0)).
\end{equation*}
$F^*(h_0)$ is a hyperbolic metric on $N$, so it's isometric to $g_0$ due to Mostow rigidity. Thus we obtain from the previous theorem the following inequality. \begin{equation*}
    \text{Area}_F(g/h_0)=d\,\text{Area}_{\text{Id}}(g/g_0)\geq d.
\end{equation*}
\end{proof}

In particular, when $n=3$, we prove a more general result for any smooth map between homeomorphic 3-manifolds which is not necessarily a diffeomorphism. Let $(M,h_0)$ and $(N,g_0)$ be closed hyperbolic 3-manifolds.
\begin{Th}\label{Thm_gromov_3}
If $\pi_1(M)\cong \pi_1(N)$, then there exists a small neighborhood $\mathcal{U}$ of $g_0$ in the $C^2$-topology, such that for any metric $g\in\mathcal{U}$ with $R_g\geq -6$, and for any smooth map $F:(N,g)\rightarrow (M,h_0)$, 
we have \begin{equation*}
    \text{deg}\,F=1 \quad \text{and}\quad \text{Area}_F(g/h_0)\geq 1.
\end{equation*}
\end{Th}

In a different situation, suppose that $\pi_1(N)$ is isomorphic to an index $d$ subgroup $G<\pi_1(M)$. Let $\Tilde{M}$ denote the covering space of $M$ with $\pi_1(\Tilde{M})=G$ and let $p:\Tilde{M}\rightarrow M$ denote the covering map, the hyperbolic metric on $\Tilde{M}$ is still represented by $h_0$. We have the following corollary. 
\begin{Cor}
There exists a small neighborhood $\mathcal{U}$ of $g_0$ in the $C^2$-topology, such that for any metric $g\in\mathcal{U}$ with $R_g\geq -6$, and for any smooth map $F:(N,g)\rightarrow (M,h_0)$ satisfying $F_*\pi_1(N)<p_*\pi_1(\Tilde{M})$, 
we have \begin{equation*}
    \text{deg}\,F=d \quad \text{and}\quad \text{Area}_F(g/h_0)\geq d.
\end{equation*}
\end{Cor}

\begin{proof}
Since $F_*\pi_1(N)<p_*\pi_1(\Tilde{M})$, $F$ can be lifted to a smooth map $\Tilde{F}:N \rightarrow \Tilde{M}$, so that $p\circ \Tilde{F}=F$. Applying Theorem \ref{Thm_gromov_3} to $\Tilde{F}$, we have \begin{equation*}
    \text{deg}\,\Tilde{F}=1 \quad \text{and}\quad \text{Area}_{\Tilde{F}}(g/h_0)\geq 1.
\end{equation*}
Thus, \begin{equation*}
    \text{deg}\,F=d \quad \text{and}\quad \text{Area}_F(g/h_0)=d\,\text{Area}_{\Tilde{F}}(g/h_0)\geq d.
\end{equation*}
\end{proof}

Besides, when the dimension of $M$ is an odd number, Hamenst{\"a}dt \cite{Hamenstadt} verified the existence of surface subgroups, and based on the result, the author discussed the corresponding minimal surface entropy and proved that $E(h_0)=2$ in \cite{Jiang}. As a result, we can extend Theorem 1.2 of \cite{LN2021} to any odd dimensions at least 3.
\begin{Th}\label{Thm_odd}
Let $(M,h_0)$ be a closed hyperbolic manifold of an odd dimension $n\geq 3$.
For any Riemannian metric $h$ on $M$, \begin{equation}
    \text{Area}_{\text{Id}}(h/h_0)E(h)\geq E(h_0)=2,
\end{equation}
the equality holds if and only if $h=ch_0$ for some constant $c>0$.
\end{Th}

This inequality can be considered as an analogue of the following comparison related to the volume entropy and geodesic stretch (Theorem 1.2 of \cite{Knieper}). More generally, suppose that $(M,h_0)$ is a compact space of any dimension with negative curvature, and $h$ is another metric on $M$. It says that \begin{equation*}
    I_{\mu_0}(h/h_0)E_{vol}(h)\geq E_{vol}(h_0),
\end{equation*}
where $\mu_0$ is the Bowen-Margulis measure of $h_0$ which is the unique measure with maximal measure theoretic entropy. The equality holds if and only if the geodesic flows of $h$ and $h_0$ are time preserving conjugate after rescaling.

We briefly mention a similar result. Besson-Courtois-Gallot \cite{Besson-Courtois-Gallot-2} showed that for any compact hyperbolic manifold (it also holds for compact locally symmetric spaces) $(M,h_0)$ with dimension $n\geq 3$, \begin{equation*}
    \big(\frac{\text{vol}_h(M)}{\text{vol}_{h_0}(M)}\big)^{\frac{1}{n}}E_{vol}(h)\geq E_{vol}(h_0),
\end{equation*}
the inequality is sharp unless $h$ is isometric to $h_0$ after rescaling. Calegari-Marques-Neves \cite{CMN} conjectured the minimal surface analogue for closed hyperbolic 3-manifolds (also see \cite{Jiang} for higher odd-dimensional closed hyperbolic manifolds and locally symmetric spaces). \begin{equation*}
    \big(\frac{\text{vol}_h(M)}{\text{vol}_{h_0}(M)}\big)^{\frac{2}{n}}E(h)\geq E(h_0)=2,
\end{equation*}
the equality holds if and only if $h$ and $h_0$ are isometric up to scaling.

The author believes that Theorem \ref{Thm_gromov_high_codim} and Theorem \ref{Thm_odd} also hold for the other compact locally symmetric spaces of rank one, which include complex hyperbolic manifolds, quaternionic hyperbolic manifolds, and the Cayley plane. To check the details, we refer the readers to the construction of the surface subgroups of locally symmetric spaces and the discussion of the equidistribution property in \cite{Kahn-Labourie-Mozes} and the corresponding definition of the minimal surface entropy in Section 5 of \cite{Jiang}.

The organization of this paper is as follows. Firstly, in Section \ref{section_preliminaries}, we introduce notations and definitions that are used throughout this paper. In Section \ref{section_gromov_high_dim}, we discuss the average area ratio for dimension $n\geq 3$ and prove Theorem \ref{Thm_gromov_high_codim}. In Section \ref{section_equidistribution}, we establish the equidistribution property and average area ratio formula which will be used in the last two sections. And in Section \ref{section_gromov_3}, we consider the special case for 3-manifolds and give the proof for Theorem \ref{Thm_gromov_3}. 
Finally, Section \ref{section_odd} concentrates on the minimal surface entropy defined on odd-dimensional manifolds, and the proof of Theorem \ref{Thm_odd} will be given.

\subsection*{Acknowledgements}
I would like to thank my advisor Andr\'{e} Neves for his constant encouragement and all the useful suggestions related to this work.

\section{Preliminaries}\label{section_preliminaries}
\subsection{Gromov's average area ratio}
Let $(M,h_0)$ be a closed hyperbolic manifold of dimension $n\geq 3$, and let $(N,g)$ be another closed Riemannian manifold of the same dimension. Suppose that $F:(N,g)\rightarrow (M,h_0)$ is a smooth map. For any $(x,L)$ in the Grassmannian bundle $Gr_2(N)$, $|\Lambda^2 F(x,L)|_g$ denotes the Jacobian of $dF_x$ at plane $L$. Take an arbitrary regular value $y\in M$ of $F$, for $(y,P)\in Gr_2(M)$, we let \begin{equation}\label{def_jac}
    |\Lambda^2 F|^{-1}_{g}(y,P) =\sum_{x\in F^{-1}(y)}\frac{1}{|\Lambda^2F(x,(dF_x)^{-1}(P)|_g}.
\end{equation}
This definition (\cite{LN2021}) is equivalent to Gromov's definition (\cite{Gromov}):
\begin{equation*}
    |\Lambda^2 F|^{-1}_{g}(y,P) =\sum_{x\in F^{-1}(y)}\lim_{\delta\rightarrow 0}\frac{\text{area}_g((dF_x)^{-1}(D_\delta))}{\delta},
\end{equation*}
where $D_\delta$ is a subset of the totally geodesic disc $D\subset \mathbb{H}^n$ which is tangential to $P$ at $x$, and $D_\delta$ has area equal to $\delta$. 

The \emph{average area ratio of $F$} is defined in \cite{Gromov} by \begin{equation}\label{def_area_ratio}
    \text{Area}_F(g/h_0)=\int_{(y,P)\in Gr_2(M)}|\Lambda^2 F|^{-1}_{g}(y,P) \,d\mu_{h_0},
\end{equation}
where $\mu_{h_0}$ stands for the unit volume measure on $Gr_2(M)$ with respect to the metric induced by $h_0$.

\subsection{Normalized Total Scalar Curvature}\label{section_E}
Let $(M,h_0)$ be a closed hyperbolic manifold of dimension $n\geq 3$, and let $\mathcal{M}$ be the space of Riemannian metrics on $M$. The \emph{total scalar curvature} (or \emph{Einstein-Hilbert functional}) $\Tilde{\mathcal{E}}: \mathcal{M}\rightarrow \mathbb{R}$ is \begin{equation*}
    \Tilde{\mathcal{E}}(h)=\int_M R_h\,dV_h.
\end{equation*}
It is a \emph{Riemannian functional} in the sense that it's invariant under diffeomorphisms, but it's not scale-invariant. To resolve this issue, we consider 
\begin{equation*}
   \mathcal{E}(h)=(\text{vol}_{h}(M))^{\frac{2}{n}-1}\int_MR_{h}\,dV_{h}.
\end{equation*}
This is called the \emph{normalized total scalar curvature} (or \emph{normalized Einstein-Hilbert functional}) of $M$. 
Under conformal deformations, the first variation of $\mathcal{E}$ is 
\begin{equation*}
    \mathcal{E}'(h)\cdot l=\frac{n-2}{2n}\text{vol}_{h}(M)^{\frac{2}{n}-1}\int_M \langle R_h-\fint_M R_h\,dV_{h},\, l\rangle_h\,dV_h.
\end{equation*}
It equals to zero provided that $(M,h)$ has constant scalar curvature. Assuming $R_h$ is constant, we can simplify the full variation to \begin{equation*}
    \mathcal{E}'(h)\cdot l=\text{vol}_{h}(M)^{\frac{2}{n}-1}\int_M \langle \frac{1}{n}R_h-Ric_h,l\rangle_h\,dV_h.
\end{equation*}
Thus, a metric $h$ is critical if and only if $(M,h)$ is Einstein. In particular, the hyperbolic metric $h_0$ is a critical point for $\mathcal{E}$.
Furthermore, since $\mathcal{E}$ is scale-invariant, from now on we may assume that 
for $h_t=h_0+tl$,  \begin{equation*}
    \int_M \frac{d}{dt}|_{t=0}(\sqrt{\text{det}_{h_0}(h_t)})\,dV_{h_0}=\frac{1}{2}\int_M \text{tr}_{h_0}l\,dV_{h_0}=0.
\end{equation*} 
Taking this into account, we obtain the second variation of $\mathcal{E}$ at $h_0$ as follows. \begin{align}\label{E_2nd_variation}
   \mathcal{E}''(h_0)(l,l)=& \text{vol}_{h_0}(M)^{\frac{2}{n}-1}\int_M \frac{d^2}{dt^2}|_{t=0}R_{h_t}-2(n-1)\frac{d^2}{dt^2}|_{t=0}(\sqrt{\text{det}_{h_0}(h_t)})\\\nonumber
    &+2\frac{d}{dt}|_{t=0}R_{h_t}\,\frac{d}{dt}|_{t=0}(\sqrt{\text{det}_{h_0}(h_t)})\,dV_{h_0},
\end{align}
Substituting the formulas 
\begin{equation*}
    \frac{d^2}{dt^2}\sqrt{\text{det}_{h_0}(h_t)}=\big(\frac{1}{4}(\text{tr}_{h_t}l)^2-\frac{1}{2}\text{tr}_{h_t}(l^2)\big)\sqrt{\text{det}_{h_0}(h_t)},
\end{equation*}
\begin{equation*}
    \frac{d}{dt}R_{h_t}=-\Delta_{h_t}(\text{tr}_{h_t}l)+\delta_{h_t}^2l-\langle Ric_{h_t}, l\rangle_{h_t},
\end{equation*}
\begin{equation*}
    \frac{d}{dt}Ric_{h_t}=-\frac{1}{2}\Delta_{h_t}l+\frac{1}{2}Ric_{h_t}(l)+\frac{1}{2}l(Ric_{h_t})-Rm_{h_t}*l-\delta^*_{h_t}(\delta_{h_t}l)-\frac{1}{2}\nabla_{h_t}^2(\text{tr}_{h_t}l),
\end{equation*}
where $(Rm_{h_t}*l)_{ij}=R_{ikjm}l^{km}$, we obtain that \begin{align*}
      \mathcal{E}''(h_0)(l,l)
    =& \text{vol}_{h_0}(M)^{\frac{2}{n}-1}\int_{M}\langle \frac{1}{2}\Delta_{h_0}l+Rm_{h_0}*l+\delta^*_{h_0}(\delta_{h_0}l)\\\nonumber
    &+\frac{n-1}{2}(\text{tr}_{h_0}l)h_0-\frac{1}{2}(\Delta_{h_0}(\text{tr}_{h_0}l))h_0+(\delta_{h_0}^2l)h_0\, ,l\rangle_{h_0} \,dV_{h_0}.
\end{align*}
According to Ebin's Slice Theorem (see \cite{Ebin}), for any $h\in\mathcal{M}$ lying in a small neighborhood of $h_0$, there exist $\phi\in\text{Diff}(M)$, $f\in C^\infty(M)$, and a \emph{transverse-traceless} tensor $l_{TT}$, i.e., $\delta_{h_0}(l_{TT})=0$ and $\text{tr}_{h_0}(l_{TT})=0$, such that $\phi^*h=h_0+fh_0+l_{TT}$. Then we can simplify the second variation.
\begin{equation}\label{E_fh_0}
    \mathcal{E}''(h_0)(fh_0,fh_0)=-\frac{(n-1)(n-2)}{2}\text{vol}_{h_0}(M)^{\frac{2}{n}-1}\int_{M}\langle \Delta_{h_0}f-nf,\, f\rangle_{h_0}\,dV_{h_0}.
\end{equation}
And we also get
\begin{equation}\label{E_l_TT}
   \mathcal{E}''(h_0)(l_{TT},l_{TT})=\text{vol}_{h_0}(M)^{\frac{2}{n}-1}\int_{M}\langle\frac{1}{2}\Delta_{h_0}l_{TT}+Rm_{h_0}*l_{TT},l_{TT}\rangle _{h_0}\,dV_{h_0}.
\end{equation}
Thus, 
\begin{align*}
    \mathcal{E}(h)&=\mathcal{E}(\phi^*h)=\mathcal{E}(h_0+fh_0+l_{TT})\\
    &=\mathcal{E}(h_0)+\mathcal{E}''(h_0)(fh_0,fh_0)+\mathcal{E}''(h_0)(l_{TT},l_{TT})+\text{higher order variations}.
\end{align*}
Taking the estimates of (\ref{E_fh_0}) and (\ref{E_l_TT}) into consideration, we can use this expansion to discuss the local behavior of $\mathcal{E}$. 
In particular, as discussed in  \cite{Besson-Courtois-Gallot}, $\mathcal{E}$ reaches a local maximum at $h_0$, and there exists $C>0$, so that for any metric $h$ in a small neighborhood of $h_0$ in $\mathcal{M}$, \begin{equation*}
    \mathcal{E}(h_0)-\mathcal{E}(h)\geq C\, d(h,h_0)^2,\quad\text{where }d(h,h_0)=\underset{\phi\in\text{Diff}(M)}{\inf} |\phi^*h-h_0|_{H^1(M,h_0)}.
\end{equation*}
The inequality is sharp unless $h$ and $h_0$ are isometric via some $\phi\in \text{Diff}(M)$.

\subsection{Equidistribution}
Suppose that the hyperbolic manifold $M$ has an odd dimension $n\geq 3$.
According to the construction by Hamenst{\"a}dt \cite{Hamenstadt}, for any small number $\epsilon>0$, there is an essential surface $\Sigma_\epsilon$ in $M$ which is sufficiently well-distributed and \emph{$(1+\epsilon)$-quasigeodesic} (i.e. the geodesics on the surface with respect to intrinsic distance are $(1+\epsilon,\epsilon)$-quasigeodesics in $M$). 
And as discussed in \cite{Jiang}, $\Sigma_\epsilon$ determines an $(1+O(\epsilon))$-quasiconformal map on $S_\infty^{n-1}$, and thus associated with an element in $S_\epsilon(M)$.

Furthermore, let $G(M,g,\epsilon)$ denote the subset of $S(M,g,\epsilon)$ consisting of homotopy classes of finite covers of $\Sigma_\epsilon$ that have genus at most $g$. 
It has cardinality comparable to $g^{2g}$. Moreover, let $S_i$ denote the minimal representative of an element in $G(M,g_i,\frac{1}{i})$, then it is homotopic to an $(1+\frac{1}{i})$-quasigeodesic surface $\Sigma_i$. From Lemma 4.3 of \cite{CMN}, for any continuous function $f$ on $M$, the unit Radon measure induced by integration over $S_i$ satisfies
\begin{equation}\label{measure}
    \underset{i\rightarrow\infty}{\lim}\,\frac{1}{\text{area}_{h_0}(S_i)}\int_{S_i} f\,dA_{h_0}=\nu(f),
\end{equation}
where the limiting measure $\nu$ is positive on any non-empty open set of $M$.

Notice that the measure $\nu$ is not necessarily identified with the unit Radon measure on $M$ induced by integration over $M$, the latter measure is denoted by $\mu$. However, in this paper, we need to find a sequence of minimal surfaces whose Radon measures defined above converge to $\mu$. To solve this problem, we introduce Labourie's construction \cite{Labourie} in Section \ref{section_equidistribution}. And we stress that both methods of Hamenst{\"a}dt and Labourie require that $M$ has an odd dimension.

\section{Gromov's Conjecture in Higher Dimensions}\label{section_gromov_high_dim}

Throughout this section, we can choose the dimension of $M$ to be any integer $n\geq 3$. The proof of Theorem \ref{Thm_gromov_high_codim} separates into two parts. If $h$ is an arbitrary metric conformal to a metric $\bar{h}$ with constant scalar curvature, we compare their average area ratios in Theorem \ref{thm_conformal}. And if $h$ is a metric with constant scalar curvature different from $h_0$, we make use of the evaluations of normalized Einstein-Hilbert functional in \cite{Besson-Courtois-Gallot}.

\subsection{Conformal deformations}

Firstly, given a metric $h$ on $M$, we look at the conformal class of $h$. Since $M$ admits a hyperbolic metric, every conformal class $[h]$ must be \emph{scalar negative}, i.e., it has a metric with negative scalar curvature. And according to the Yamabe problem, after rescaling, there exists a unique metric $\bar{h}\in [h]$ with constant scalar curvature $c<0$. 
In Theorem \ref{Thm_gromov_high_codim}, we assume that $c\geq -n(n-1)$. Set $\bar{h}'=\frac{c}{-n(n-1)}\bar{h}$, so $R_{\bar{h}'}\equiv -n(n-1)$, and for any surface or 2-dimensional subset $S$ in $M$, we have 
\begin{equation*}
    \text{area}_{\bar{h}'}(S)=\frac{c}{-n(n-1)}\text{area}_{\bar{h}}(S)\leq\text{area}_{\bar{h}}(S).
\end{equation*}
Therefore, to prove the theorem, we may assume that $R_{\bar{h}}\equiv -n(n-1)$.

\begin{Th}\label{thm_conformal}
Suppose the scalar curvature $R_h\geq -n(n-1)$, then we have \begin{equation}\label{conformal_area_ratio}
    \text{Area}_{\text{Id}}(h/h_0)\geq \text{Area}_{\text{Id}}(\bar{h}/h_0).
\end{equation}
Furthermore, 
for any surface subgroup $\Pi\in S_\epsilon(M)$,\begin{equation}\label{conformal_minimal_area}
    \text{area}_h(\Pi)\geq\text{area}_{\bar{h}}(\Pi),
\end{equation}
and as an immediate result,\begin{equation}\label{conformal_entropy}
    E(h)\leq E(\bar{h}).
\end{equation}
Each of the above equalities holds if and only if $h=\bar{h}$.
\end{Th}

\begin{proof}

Set $h=e^{2\phi}\bar{h}$. The conformal factor $\phi$ satisfies that \begin{align*}
    e^{2\phi} R_h&=R_{\bar{h}}-2(n-1)\Delta_{\bar{h}}\phi-(n-2)(n-1)|d\phi|_{\bar{h}}^2\\
    &=-n(n-1)-2(n-1)\Delta_{\bar{h}}\phi-(n-2)(n-1)|d\phi|_{\bar{h}}^2.
\end{align*}
Let $\phi_{min}:=\underset{x\in M}{\min}\,\phi(x)$, we obtain $\Delta_{\bar{h}}\phi_{min}\geq 0$, and it yields that \begin{equation}\label{conformal_inequ}
    -n(n-1)e^{2\phi_{min}}\leq e^{2\phi_{min}}R_h\leq -n(n-1).
\end{equation}
Thus, for any subset $D_\delta$ of a totally geodesic disc in $\mathbb{H}^n$ with hyperbolic area equal to $\delta$, we have \begin{equation*}
    e^{2\phi_{min}}\geq 1\quad \Longrightarrow \quad \text{area}_h(D_\delta)=\int_{D_\delta}e^{2\phi}\,dA_{\bar{h}} \geq \text{area}_{\bar{h}}(D_\delta).
\end{equation*}
The inequality (\ref{conformal_area_ratio}) follows from the definition of the average area ratio (\ref{def_area_ratio}). 

In addition, it also shows that for any surface $S\in M$, \begin{equation*}
    \text{area}_h(S)=\int_{S}e^{2\phi}\,dA_{\bar{h}} \geq \text{area}_{\bar{h}}(S).
\end{equation*}
Thus we conclude (\ref{conformal_minimal_area}), and the comparison of entropy (\ref{conformal_entropy}) is a direct corollary of the definition (\ref{eh}).

Moreover, if any equality in the theorem holds, the inequality \ref{conformal_inequ} implies that $R_h\equiv -n(n-1)$, due to the uniqueness of the solution to Yamabe problem, $h$ must be identical to $\bar{h}$.
\end{proof}

\subsection{Definition of functional $\mathcal{A}$}

Let $\mathcal{M}$ be the space of all Riemannian metrics on $M$, and let $\mathcal{M}_R$ be the subset consisting of metrics with constant scalar curvature $-n(n-1)$. 
From the previous section, it remains to consider metrics in $\mathcal{M}_R$. 

Define a functional $\mathcal{A}$ from the space of Riemannian metrics on $M$ to $\mathbb{R}$ as follows. \begin{align*}
    \mathcal{A}(h)=& \int_{(x,P)\in Gr_2M}R_{h}(x)\,|\Lambda^2\text{Id}|^{-1}_h(x,P) \,d\mu_{h_0}\\
    =& \lim_{\delta\rightarrow 0}\int_{(x,P)\in Gr_2M} R_{h}(x)\,\frac{\text{area}_h(D_\delta(P))}{\delta}\,d\mu_{h_0}\\
    =& \lim_{\delta\rightarrow 0}\int_{(x,P)\in Gr_2M} R_{h}(x)\,\fint_{D_\delta(P)}\sqrt{\text{det}_{h_0}(h|_{D_\delta(P)})}\,dA_{h_0}\,d\mu_{h_0},
\end{align*}
where $D_\delta(P)$ is a subset of the totally geodesic disc in $\mathbb{H}^n$ that is tangential to $P$ at $x$, and $D_\delta(P)$ has hyperbolic area equal to $\delta$, $\mu_{h_0}$ is the unit volume measure on $Gr_2M$ with respect to the metric induced by $h_0$.
Notice that the metric $h\in\mathcal{M}_R$ satisfies  \begin{align*}
    \mathcal{A}(h)&=-n(n-1)\,\int_{(x,P)\in Gr_2M}|\Lambda^2\text{Id}|^{-1}_h(x,P) \,d\mu_{h_0}\\ 
    &=\mathcal{A}(h_0)\,\int_{(x,P)\in Gr_2M}|\Lambda^2\text{Id}|^{-1}_h(x,P) \,d\mu_{h_0}.
\end{align*}
Therefore, from the definition of the average area ratio (\ref{def_area_ratio}), to deduce that if $h$ is in a small neighborhood of $h_0$ in $\mathcal{M}_R$, then \begin{equation*}
  \text{Area}_{\text{Id}}(h/h_0)=\int_{(x,P)\in Gr_2M}|\Lambda^2\text{Id}|^{-1}_h(x,P) \,d\mu_{h_0}\geq 1,
\end{equation*}
we only need to show that $\mathcal{A}$ attains a local maximum at $h_0$. 


\subsection{Proof of Theorem \ref{Thm_gromov_high_codim}}
To see this, we'll discuss the first and second variations of $\mathcal{A}$ in detail.  Before start, we notice that the functional $\mathcal{A}$ is scale-invariant, so it suffices to assume $\text{vol}_{h_0}(M)=1$, and the symmetric tensor $l=h-h_0$ satisfies $\int_M \text{tr}_{h_0}l\,dV_{h_0}=0$. Let $h_t=h_0+tl$,
we rewrite $\mathcal{A}(h_t)$ as \begin{equation*}
    \mathcal{A}(h_t)=\int_{x\in M} R_{h_t}(x) \,a_{h_t}(x)\,dV_{h_0},
\end{equation*} 
where \begin{equation*}
    a_{h_t}(x)=\lim_{\delta\rightarrow 0}\fint_{P\in Gr_2M_x} \fint_{D_\delta(P)}\sqrt{\text{det}_{h_0}(h_t|_{D_\delta(P)})}\,dV_{h_0}d\nu_{h_0}.
\end{equation*}
Then we have 
\begin{equation}\label{1st_variation}
    \frac{d}{dt}\mathcal{A}(h_t)
    =\int_M \frac{d}{dt}R_{h_t}\,a_{h_t}+R_{h_t}\frac{d}{dt}a_{h_t}\,dV_{h_0}.
\end{equation}
In addition, \begin{equation}\label{d_R}
    \frac{d}{dt}R_{h_t}=-\Delta_{h_t}(\text{tr}_{h_t}l)+\delta_{h_t}^2l-\langle Ric_{h_t}, l\rangle_{h_t},
\end{equation}
where $\Delta_{h_t}$ is the rough Laplacian with negative eigenvalues, and $\delta_{h_t}^2$ is the double divergence operator.
In particular, when $t=0$, applying the Stokes' theorem, we have \begin{equation*}
    \int_M\frac{d}{dt}|_{t=0}\,R_{h_t}\,a_{h_0}\,dV_{h_0}= \int_M -\langle Ric_{h_0},l\rangle_{h_0}\,dV_{h_0}=(n-1)\int_M\text{tr}_{h_0}l\,dV_{h_0}.
\end{equation*}
On the other hand, since
\begin{equation*}
  \frac{d}{dt}\sqrt{\text{det}_{h_0}(h_t|_{D_\delta(P)})} =\frac{1}{2}\text{tr}_{h_t}l|_{D_\delta(P)} \,\sqrt{\text{det}_{h_0}(h_t|_{D_\delta(P)})},
\end{equation*}
let $\lambda_1,\cdots,\lambda_n$ be the eigenvalues of the matrix $l$ at point $x\in M$ with respect to $h_0$, we obtain by computation that \begin{align*}
    \frac{d}{dt}|_{t=0}\,a_{h_t} &=\lim_{\delta\rightarrow 0}\frac{1}{2}\fint_{P\in Gr_2M_x}\fint_{D_\delta(P)}\text{tr}_{h_0}l|_{D_\delta(P)}\,dA_{h_0}d\nu_{h_0}\\ 
    &= \frac{1}{2}\,\frac{\sum_{i< j}\lambda_i+\lambda_j}{\binom{n}{2}} = \frac{1}{2}\,\frac{(n-1)\sum_{i=1}^n\lambda_i}{\binom{n}{2}}\\ 
    &=\frac{1}{n}\text{tr}_{h_0}l.
\end{align*}
It follows that for any symmetric 2-tensor $l$,
\begin{equation*}
    \mathcal{A}'(h_0)\cdot l =(n-1)\int_M\text{tr}_{h_0}l\,dV_{h_0}-n(n-1)\frac{1}{n}\int_M \text{tr}_{h_0}l\,dV_{h_0}
    =0,
\end{equation*}
thus $h_0$ is a critical point of $\mathcal{A}$.

Now we proceed to compute the second variation at $t=0$. Note that $\mathcal{A}$ is an analogue of the normalized total scalar curvature $\mathcal{E}$, it's easier to compare their second variations using the computation in Section \ref{section_E}.
Based on (\ref{1st_variation}) and (\ref{E_2nd_variation}),
\begin{align}\label{2nd_variation_1}
    &\mathcal{A}''(h_0)(l,l)-\mathcal{E}''(h_0)(l,l) \\\nonumber
    =& \int_M -n(n-1)\frac{d^2}{dt^2}|_{t=0}\,a_{h_t}+2(n-1)\frac{d^2}{dt^2}|_{t=0}(\sqrt{\text{det}_{h_0}(h_t)})\,dV_{h_0}\\\nonumber
    &+ \int_M 2\frac{d}{dt}|_{t=0}R_{h_t}(\frac{d}{dt}|_{t=0}(a_{h_t}-\sqrt{\text{det}_{h_0}(h_t)})\,dV_{h_0}\\\nonumber
    =& \int_M -n(n-1)\frac{d^2}{dt^2}|_{t=0}\,a_{h_t}+2(n-1)\big(\frac{1}{4}(\text{tr}_{h_0}l)^2-\frac{1}{2}\text{tr}_{h_0}(l^2)\big)\,dV_{h_0}\\\nonumber
    &+2\int_M \big(-\Delta_{h_0}(\text{tr}_{h_0}l)+\delta_{h_0}^2l+(n-1)\text{tr}_{h_0}l\big)(\frac{1}{n}-\frac{1}{2})\text{tr}_{h_0}l\,dV_{h_0}.
\end{align}

Next, using \begin{equation*}
    \frac{d^2}{dt^2}\sqrt{\text{det}_{h_0}(h_t|_{D_\delta(P)})}=\big(\frac{1}{4}(\text{tr}_{h_t}l|_{D_\delta(P)})^2-\frac{1}{2}\text{tr}_{h_t}(l|_{D_\delta(P)}^2)\big)\sqrt{\text{det}_{h_0}(h_t|_{D_\delta(P)})},
\end{equation*}
we estimate the first term on the right-hand side of (\ref{2nd_variation_1}). \begin{align}\label{limit_estimate}
&-n(n-1)\int_M\frac{d^2}{dt^2}|_{t=0}\,a_{h_t}\,dV_{h_0}\\\nonumber
    =&-n(n-1)\lim_{\delta\rightarrow 0}\,\int_{x\in M} \fint_{P\in Gr_2M_x} \fint_{D_\delta(P)}\frac{1}{4}(\text{tr}_{h_0}l|_{D_\delta(P)})^2-\frac{1}{2}\text{tr}_{h_0}(l|_{D_\delta(P)}^2)\,dA_{h_0}d\nu_{h_0}dV_{h_0}\\ \nonumber
    =& -n(n-1)\int_M \frac{\sum_{i< j}\frac{1}{4}(\lambda_i+\lambda_j)^2-\frac{1}{2}(\lambda_i^2+\lambda_j^2)}{\binom{n}{2}}\,dV_{h_0}\\ \nonumber
    =& -n(n-1)\int_M \frac{\frac{1}{4}(\sum_{i=1}^n\lambda_i)^2-\frac{n}{4}\sum_{i=1}^n\lambda_i^2}{\binom{n}{2}}\,dV_{h_0}\\ \nonumber
    =& \int_M -\frac{1}{2}(\text{tr}_{h_0}l)^2+\frac{n}{2}\text{tr}_{h_0}(l^2)\,dV_{h_0}.
\end{align}
Combining (\ref{2nd_variation_1}) and (\ref{limit_estimate}), we have \begin{align}\label{A_2nd_variation}
     \mathcal{A}''(h_0)(l,l) =& \mathcal{E}''(h_0)(l,l)- \int_M \frac{n-2}{2}\text{tr}_{h_0}(l^2) \,dV_{h_0}\\\nonumber
     &-\int_M\frac{(n-2)^2}{2n}(\text{tr}_{h_0}l)^2
     +\frac{n-2}{n}|\nabla_{h_0}(\text{tr}_{h_0}l)|^2+\frac{n-2}{n}\delta_{h_0}^2l\,\text{tr}_{h_0}l\,dV_{h_0}.
\end{align}

To simplify this quadratic form, we decompose $l$ into three parts.
Applying the decomposition of space of symmetric tensors for a compact Einstein manifold other than the standard sphere (Theorem 4.60 of \cite{Besse}), we have \begin{equation*}
    T_{h_0}\mathcal{M}=C^\infty(M)\cdot h_0\oplus T_{h_0}(\text{Diff}(M)(h_0))\oplus TT_{h_0}, 
\end{equation*}
where $C^\infty(M)\cdot h_0$ represents the conformal deformations of $h_0$, $\text{Diff}(M)(h_0)$ is the action of the diffeomorphism group on $h_0$, and $TT_{h_0}=\ker \text{tr}_{h_0}\cap \ker \delta_{h_0}$ stands for the set of transverse-traceless tensors. 
Let $h\in\mathcal{M}_R$, and let $l=h-h_0$, it decomposes into $l=fh_0+l_D+l_{TT}$, where $fh_0\in C^\infty(M)\cdot h_0$, $l_D\in T_{h_0}(\text{Diff}(M)(h_0))$, and $l_{TT}\in TT_{h_0}$.
The second variation of $\mathcal{A}$ has the form
\begin{align*}
  &\mathcal{A}''(h_0)(l,l)\\
   =&\mathcal{A}''(h_0)(l_D+l_{TT}, l_D+l_{TT})+\mathcal{A}''(h_0)(fh_0,l_D)+\mathcal{A}''(h_0)(l_D,fh_0)\\
    & +\mathcal{A}''(h_0)(fh_0,l_{TT})+\mathcal{A}''(h_0)(l_{TT},fh_0)+\mathcal{A}''(h_0)(fh_0,fh_0)\\
   \leq &\mathcal{A}''(h_0)(l_D+l_{TT}, l_D+l_{TT})+C_0|fh_0|_{H^1(M,h_0)}|l|_{H^1(M,h_0)}.
\end{align*}
By Theorem \ref{thm_conformal}, to check that $\mathcal{A}$ reaches a local maximum at $h_0$ on $\mathcal{M}_R$, it remains to analyze the sign of $\mathcal{A}''(h_0)(l_D+l_{TT},l_D+l_{TT})$ and estimate $|fh_0|_{H^1(M,h_0)}$. To see these, we prove the following lemmas.

\begin{Lemma}\label{lem_3.2}
There exists a constant $C>0$, such that in the decomposition of $l$, $l_D\in T_{h_0}(\text{Diff}(M)(h_0))$ and $l_{TT}\in TT_{h_0}$ satisfy 
\begin{equation*}
    \mathcal{A}''(h_0)((l_D+l_{TT}, l_D+l_{TT})\leq -C(|l_{TT}|^2_{H^1(M,h_0)}+|X|^2_{L^2(M,h_0)}),
\end{equation*}
where $\mathcal{L}_Xh_0=l_D$. 
\end{Lemma}

\begin{Lemma}\label{lem_3.3}
For any $\epsilon>0$, we can find a $C^2$-neighborhood $\mathcal{U}_{\epsilon,R}$ of $h_0$ on $\mathcal{M}_R$ and $c>0$, such that for any $h\in\mathcal{U}_{\epsilon,R}$, 
    $fh_0\in C^\infty(M)\cdot h_0$ in the decomposition of $l=h-h_0$ satisfies \begin{equation*}
    |fh_0|_{H^1(M,h_0)}\leq c\,\epsilon|l|_{H^1(M,h_0)}.
    \end{equation*}
\end{Lemma}

\begin{proof}[Proof of Lemma \ref{lem_3.2}]
\begin{align}\label{3.2.1}
    &\mathcal{A}''(h_0)((l_D+l_{TT}, l_D+l_{TT})\\\nonumber
    =&\mathcal{A}''(h_0)(l_{TT}, l_{TT})+\mathcal{A}''(h_0)((l_D, l_D)+\mathcal{A}''(h_0)(l_D,l_{TT})+\mathcal{A}''(h_0)(l_{TT},l_D).
\end{align}
To find an upper bound of the first term, we use the estimate in Lemma 2.9 of \cite{Besson-Courtois-Gallot}. There exists a constant $C_1>0$, such that 
\begin{equation*}
   \mathcal{E}''(h_0)(l_{TT},l_{TT})=\int_{M}\langle\frac{1}{2}\Delta_{h_0}l_{TT}+Rm_{h_0}*l_{TT},l_{TT}\rangle _{h_0}\,dV_{h_0}\leq -C_1|l_{TT}|^2_{H^1(M,h_0)}.
\end{equation*}
And thus, substituting the above inequality and $\text{tr}_{h_0}l_{TT}=\delta_{h_0}l_{TT}=0$ into (\ref{A_2nd_variation}), we obtain \begin{align}\label{3.2.2}
    \mathcal{A}''(h_0)(l_{TT}, l_{TT})&=  \mathcal{E}''(h_0)(l_{TT},l_{TT})-\frac{n-2}{2}\int_M \text{tr}_{h_0}(l_{TT}^2)\,dV_{h_0}\\\nonumber
    &\leq -C_1|l_{TT}|^2_{H^1(M,h_0)}.
\end{align}

To deal with the remaining terms of (\ref{3.2.1}), we apply the diffeomorphism invariance property of $\mathcal{E}$, which says \begin{equation}\label{diffeo_invariance}
    \mathcal{E}''(h_0)(l_D,\cdot)=\mathcal{E}''(h_0)(\cdot,l_D)=0.
\end{equation}
Moreover, for any $\phi\in \text{Diff}(M)$, \begin{equation}\label{dR=0}
    R_{\phi^*h_0}\equiv R_{h_0}\equiv -n(n-1)\quad\Longrightarrow \quad R_{h_0}'\cdot l_D=0.
\end{equation}
Therefore, the second variations comparison (\ref{2nd_variation_1}), in company with (\ref{limit_estimate}), says that 
\begin{equation}\label{l_d}
    \mathcal{A}''(h_0)((l_D, l_D) 
    =\frac{n-2}{2}\int_M (\text{tr}_{h_0}l_D)^2-\text{tr}_{h_0}(l_D^2)\,dV_{h_0}.
\end{equation}
Since $l_D\in T_{h_0}(\text{Diff}(M)(h_0))$, it can be expressed by the Lie derivative of the metric $h_0$ in the direction $X$.
And using Helmholtz-Hodge decomposition, $X$ decomposes further into $X=\nabla_{h_0}r+Y$, where $r$ is a scalar function, and $Y$ is a vector field with $\text{div}_{h_0}Y=0$. By computation, \begin{equation}\label{l_d_0}
    \text{tr}_{h_0}l_D=2\Delta_{h_0}r+2\text{div}_{h_0}Y=2\Delta_{h_0}r.
\end{equation}
In addition, \begin{align}\label{l_d_1}
|l_D|^2_{L^2(M,h_0)}=&
|2\nabla^2_{h_0}r+\sum_{i,j}(\nabla_iY_j+\nabla_jY_i)|^2_{L^2(M,h_0)}\\\nonumber
=& 4|\nabla^2_{h_0}r|^2_{L^2(M,h_0)}+|\sum_{i,j}(\nabla_iY_j+\nabla_jY_i)|^2_{L^2(M,h_0)}\\\nonumber
&+4\langle\nabla^2_{h_0}r,\sum_{i,j}(\nabla_iY_j+\nabla_jY_i)\rangle_{L^2(M,h_0)}.
\end{align}
where the second term is in the form \begin{align}\label{l_d_2}
    &|\sum_{i,j}(\nabla_iY_j+\nabla_jY_i)|^2_{L^2(M,h_0)} \\\nonumber
    =& \int_M 4\sum_i (\nabla_i Y_i)^2+\sum_{i\neq j}(\nabla_iY_j)^2+(\nabla_jY_i)^2+2\nabla_i Y_j \nabla_j Y_i\,dV_{h_0}\\\nonumber
    =&\int_M 4(\text{div}_{h_0}Y)^2+\sum_{i\neq j}-2\nabla_i Y_i\nabla_jY_j+(\nabla_iY_j)^2+(\nabla_jY_i)^2\\\nonumber
    &-2\nabla_i Y_i\nabla_jY_j+2\nabla_i Y_j \nabla_j Y_i\,dV_{h_0}\\\nonumber
    =& \int_M \sum_{i\neq j} Y_j\nabla_j\nabla_iY_i+Y_i\nabla_i\nabla_jY_j+(\nabla_iY_j)^2+(\nabla_jY_i)^2\\\nonumber
    &+Y_i\nabla_i\nabla_j Y_j-Y_i\nabla_j\nabla_i Y_j +Y_j\nabla_j\nabla_i Y_i-Y_j\nabla_i\nabla_j Y_i\,dV_{h_0}\\\nonumber
    =& \int_M \sum_{i\neq j} Y_j(\nabla_i\nabla_jY_i-R_{ijj}^iY_j)+Y_i(\nabla_j\nabla_iY_j-R_{jii}^jY_i)+(\nabla_iY_j)^2+(\nabla_jY_i)^2\\\nonumber
    &+Y_i\nabla_i\nabla_j Y_j-Y_i(\nabla_i\nabla_j Y_j+R_{jii}^j Y_i) +Y_j\nabla_j\nabla_i Y_i-Y_j(\nabla_j\nabla_i Y_i+R_{ijj}^iY_j)\,dV_{h_0}\\\nonumber
    =& \int_M \sum_{i\neq j} Y_i^2+Y_j^2-2\nabla_jY_i\nabla_iY_j+(\nabla_iY_j)^2+(\nabla_jY_i)^2+Y_i^2+Y_j^2\,dV_{h_0}\\\nonumber
    =& 4(n-1)|Y|^2_{L^2(M,h_0)}+|\sum_{i, j}(\nabla_iY_j-\nabla_jY_i)|^2_{L^2(M,h_0)},
\end{align}
and the last term of (\ref{l_d_1}) vanishes, since \begin{align}\label{l_d_3}
    &4\langle\nabla^2_{h_0}r,\sum_{i,j}(\nabla_iY_j+\nabla_jY_i)\rangle_{L^2(M,h_0)}\\\nonumber
    =& 4\int_M 2\sum_i\nabla_i\nabla_ir\nabla_iY_i +\sum_{i\neq j}\nabla_i\nabla_j r(\nabla_iY_j+\nabla_j Y_i)\,dV_{h_0}\\\nonumber
    =& 4\int_M 2\sum_i\nabla_i\nabla_ir\nabla_iY_i -\sum_{i\neq j}(\nabla_ir \nabla_j\nabla_i Y_j+\nabla_j r\nabla_i\nabla_j Y_i)\,dV_{h_0}\\\nonumber
    =& 4\int_M 2\sum_i\nabla_i\nabla_ir\nabla_iY_i -\sum_{i\neq j}\big(\nabla_ir(\nabla_i\nabla_jY_j+R_{jii}^jY_i)+\nabla_jr(\nabla_j\nabla_iY_i+R_{ijj}^i Y_j)\big)\,dV_{h_0}\\\nonumber
    =&4\int_M 2\Delta_{h_0}r\,\text{div}_{h_0}Y -\sum_{i\neq j} (\nabla_i\nabla_i r\nabla_j Y_j+\nabla_j\nabla_jr\nabla_iY_i)\\\nonumber
    &+ \sum_{i\neq j} (\nabla_i\nabla_i r\nabla_j Y_j+\nabla_j\nabla_jr\nabla_iY_i) +\sum_{i\neq j} \nabla_ir Y_i+\nabla_jrY_j\,dV_{h_0}\\\nonumber
    =& 8(n-1)\langle \nabla_{h_0}r,Y\rangle_{L^2(M,h_0)}=0.
\end{align}
Substituting (\ref{l_d_0})-(\ref{l_d_3}) into (\ref{l_d}), then applying Bochner's formula, we obtain \begin{align}\label{3.2.3}
    &\mathcal{A}''(h_0)(l_D,l_D) \\\nonumber
    =& \frac{n-2}{2}\big(\int_M 4(\Delta_{h_0}r)^2-4|\nabla_{h_0}^2r|^2\,dV_{h_0}\\\nonumber
    &-4(n-1)|Y|^2_{L^2(M,h_0)}-|\sum_{i, j}(\nabla_iY_j-\nabla_jY_i)|^2_{L^2(M,h_0)}\big)\\\nonumber
    =& \frac{n-2}{2}\big(\int_M 4Ric(\nabla_{h_0}r,\nabla_{h_0}r)\,dV_{h_0}-4(n-1)|Y|^2_{L^2(M,h_0)}-|\sum_{i, j}(\nabla_iY_j-\nabla_jY_i)|^2_{L^2(M,h_0)}\big)\\\nonumber
    \leq & -2(n-1)(n-2)|X|^2_{L^2(M,h_0)}.
\end{align}
Furthermore, after eliminating some terms using (\ref{diffeo_invariance}), (\ref{dR=0}) and $tr_{h_0}l_{TT}=0$, we get \begin{equation}\label{3.2.4}
    \mathcal{A}''(h_0)(l_D,l_{TT})=-\frac{n-2}{2}\int_M \langle l_D,l_{TT}\rangle_{h_0}\,dV_{h_0}=0,
\end{equation}
it vanishes because of the $L^2$-orthogonality between $T_{h_0}(\text{Diff}(M)(h_0))$ and $TT_{h_0}$. 
Similarly, the traceless-transverse property of $l_{TT}$ simplifies (\ref{A_2nd_variation} to \begin{equation}\label{3.2.5}
    \mathcal{A}''(h_0)(l_{TT},l_D)=-\frac{n-2}{2}\int_M \langle l_{TT},l_{D}\rangle _{h_0}\,dV_{h_0}=0.
\end{equation}

Substituting (\ref{3.2.2}), (\ref{3.2.3}), (\ref{3.2.4}) and (\ref{3.2.5}) into (\ref{3.2.1}), we complete the proof. 

\end{proof}

\smallskip
\begin{proof}[Proof of Lemma \ref{lem_3.3}]

Since $h\in\mathcal{M}_R$, we have 
\begin{equation*}
        0=R_h-R_{h_0}=\int_0^1 R_{h_0+tl}'\cdot l\,dt,
    \end{equation*}
it leads to 
    \begin{equation}\label{3.3.1}
        \langle R_{h_0}'\cdot l,\text{tr}_{h_0}(fh_0)\rangle_{L^2(M,h_0)} =-\int_0^1\langle (R'_{h_0+tl}-R'_{h_0})\cdot l,\text{tr}_{h_0}(fh_0)\rangle_{L^2(M,h_0)}\,dt.
    \end{equation}
On the left-hand side, we have $R_{h_0}'\cdot l=R_{h_0}'\cdot(fh_0)$ as analyzed earlier, and it follows from (\ref{d_R}) that \begin{align}\label{3.3.2}
   & \langle R_{h_0}'\cdot l,\text{tr}_{h_0}(fh_0)\rangle_{h_0} \\\nonumber
   =& \int_M\langle -\Delta_{h_0}(\text{tr}_{h_0}(fh_0))+\delta^2_{h_0}(fh_0)+(n-1)\text{tr}_{h_0}(fh_0),\text{tr}_{h_0}(fh_0)\rangle_{h_0}\,dV_{h_0}\\\nonumber
    =& \int_M n(n-1)|\nabla_{h_0}f|^2+n^2(n-1)|f|^2\,dV_{h_0}\\\nonumber
    \geq & (n-1)|fh_0|^2_{H^1(M,h_0)}.
\end{align}
On the other hand, the right-hand side of (\ref{3.3.1}) can be estimated using the continuity of $h\rightarrow\text{tr}_{h}l$, $h\rightarrow \langle\cdot,\cdot\rangle_{L^2(M,h)}$, $h\rightarrow \nabla_h l$, and $h\rightarrow \delta_h l$. For any $\epsilon>0$, after shrinking the neighborhood of $h_0$ in $\mathcal{M}_R$, we have \begin{equation}\label{3.3.3}
-\int_0^1\langle (R'_{h_0+tl}-R'_{h_0})\cdot l,\text{tr}_{h_0}(fh_0)\rangle_{L^2(M,h_0)}\,dt  \leq \epsilon |h-h_0|_{C^2}|l|_{H^1(M,h_0)}|fh_0|_{H^1(M,h_0)}.
\end{equation}
Therefore, taken (\ref{3.3.2}) and (\ref{3.3.3}) into account, (\ref{3.3.1}) leads to the result.

\end{proof}

To end this section, we discuss the equality condition of Theorem \ref{Thm_gromov_high_codim}. There exists $t\in(0,1)$, such that \begin{equation*}
    \mathcal{A}(h)=\mathcal{A}(h_0)+\frac{\mathcal{A}(h_0)''(l,l)}{2}+\frac{\mathcal{A}'''(h_0+tl)(l,l,l)}{6}.
\end{equation*}
Following the same procedure to compute $\mathcal{A}'''(h_0)$ and discuss the continuity near $h_0$, we can see that $\mathcal{A}'''(h_0+tl)(l,l,l)=O(|l|^3_{H^1(M,h_0)})$. 
Thus, from the above lemmas, there exists $C'>0$,
\begin{equation*}
    \mathcal{A}(h)\leq \mathcal{A}(h_0) -C'(|fh_0+l_{TT}|^2_{H^1(M,h_0)}+|X|^2_{L^2(M,h_0)})+O(|h-h_0|^3_{H^1(M,h_0)}).
\end{equation*} 
From the computation of (\ref{3.2.3}) and its derivative, we can see that as $h\rightarrow h_0$, $|X|_{L^2(M,h_0)}$ and $|l_D|_{H^1(M,h_0)}$ decay at the same rate, 
so the norm $|fh_0+l_{TT}|^2_{H^1(M,h_0)}+|X|^2_{L^2(M,h_0)}\sim C''|h-h_0|^2_{H^1(M,h_0)}$.
Consequently, the following expansion holds for metric $h$ in a small neighborhood of $h_0$ with $R_h\geq -n(n-1)$. \begin{equation*}
    \text{Area}_{\text{Id}}(h/h_0) \geq  1+C'''(|fh_0+l_{TT}|^2_{H^1(M,h_0)}+|X|^2_{L^2(M,h_0)})+O(|h-h_0|^3_{H^1(M,h_0)}).
\end{equation*}
where $C'''>0$. Note that $\text{Area}_{\text{Id}}(h/h_0)=1$ requires that $fh_0=l_{TT}=X=0$, and thus $h=h_0$.

\section{Equidistribution property and average area ratio formula}\label{section_equidistribution}
In this section, we extend the notations of laminations and associated properties for hyperbolic 3-manifolds (see Labourie \cite{Labourie} and Lowe-Neves \cite{LN2021}) to the higher, odd-dimensional case. The purpose of this section is to deduce the average area ratio formula in Lemma \ref{lemma_area ratio formula}, as an important tool in the proofs of Theorem \ref{Thm_gromov_3} and Theorem \ref{Thm_odd}.

\subsection{Laminations and laminar measures}
Let $M=\mathbb{H}^n/\pi_1(M)$ be a closed hyperbolic manifold of dimension $n\geq 3$.
And let $\mathcal{F}(\mathbb{H}^n,\epsilon)$ be the space of conformal minimal immersions $\Phi:\mathbb{H}^2\rightarrow\mathbb{H}^n$,
such that $\Phi(\partial_\infty\mathbb{H}^2)$ is an $(1+\epsilon)$-quasicircle. As discussed in Section 2 of \cite{Jiang}, when $\epsilon$ is sufficiently small, $\Phi(\mathbb{H}^2)$ is a stable embedded disc in $\mathbb{H}^n$. The space $\mathcal{F}(\mathbb{H}^n,\epsilon)$ equips with the topology of uniform convergence on compact sets, and we take \begin{equation*}
    \mathcal{F}(M,\epsilon):=\mathcal{F}(\mathbb{H}^n,\epsilon)/\pi_1(M)
\end{equation*} with the quotient topology. The space $\mathcal{F}(M,\epsilon)$ together with the action of $\text{PSL}(2,\mathbb{R})$ by pre-composition \begin{equation}\label{psl2r}
    \mathcal{R}_\gamma:\mathcal{F}(M,\epsilon)\rightarrow\mathcal{F}(M,\epsilon),\quad \mathcal{R}_\gamma(\phi)=\phi\circ\gamma^{-1},\quad\forall\gamma\in \text{PSL}(2,\mathbb{R})
\end{equation}
is called the \emph{conformal minimal lamination of $M$}. A \emph{laminar measure} on $\mathcal{F}(M,\epsilon)$ stands for a probability measure which is invariant under the $\text{PSL}(2,\mathbb{R})$-action defined as above. 
The space $\mathcal{F}(M,\epsilon)$ is sequentially compact, but the space of laminar measures is not necessarily weakly compact. In light of that, we consider a continuous map from $\mathcal{F}(M,\epsilon)$ to the frame bundle $F(M)$ of $M$, the latter space is compact, so the space of probability measures on $F(M)$ is compact in weak-$*$ topology. 

Firstly we define a map from $\mathcal{F}(M,\epsilon)$ to the 2-vector bundle $F_2(M)$ on $M$ consisting of $(x,v_1,v_2)\in M\times S_xM\times S_xM$, where $S_xM$ denotes the unit sphere in the tangent space to $M$ at $x$. Let $\{e_1,e_2\}$ be an orthonormal basis of $\mathbb{H}^2$, and for any $\phi\in\mathcal{F}(M,\epsilon)$, let $\phi^*(h_0)=C_\phi^2h_{\mathbb{H}^2}$, where $C_\phi^2>0$ denotes the conformal factor between the hyperbolic metric on $\mathbb{H}^2 $ and the pull-back metric of $h_0$ by $\phi$. Let $e_1(\phi)=\dfrac{d\phi(e_1)}{C_\phi}$ and $e_2(\phi)=\dfrac{d\phi(e_2)}{C_\phi}$. We define the following continuous map.
\begin{equation*}
    \Omega:\mathcal{F}(M,\epsilon)\rightarrow F_2(M),\quad \Omega(\phi)=(\phi(i),e_1(\phi),e_2(\phi)).
\end{equation*}
Furthermore, it induces a map from $\mathcal{F}(M,\epsilon)$ to the frame bundle $F(M)$ by parallel transport: \begin{equation*}
     \overline{\Omega}:\mathcal{F}(M,\epsilon)\rightarrow F(M),\quad \overline{\Omega}(\phi)=(\phi(i),e_1(\phi),e_2(\phi),\cdots,e_n(\phi)).
\end{equation*}
And define the projection \begin{equation*}
    P:F(M)\rightarrow F_2(M),\quad P(x,e_1,e_2,\cdots,e_n)=(x,e_1,e_2).
\end{equation*}

We consider the subspace $\mathcal{F}(\mathbb{H}^n,0)\subset \mathcal{F}(\mathbb{H}^n,\epsilon)$, it contains isometric immersions $\phi_0:\mathbb{H}^2\rightarrow \mathbb{H}^n$ whose images are totally geodesic discs in $\mathbb{H}^n$. Conversely, each totally geodesic disc is uniquely determined by $\phi(i)$, and tangent vectors $e_1(\phi)$, $e_2(\phi)$. 
Let $\Omega_0:\mathcal{F}(M,0)\rightarrow F_2(M)$ be the restriction of $\Omega$ to $\mathcal{F}(M,0)$, it's therefore a bijection. Using (\ref{psl2r}), We can define the $\text{PSL}(2,\mathbb{R})$-action on $F(M)$ as follows. \begin{equation*}
     R_\gamma:F(M)\rightarrow F(M),\quad R_\gamma(x)=\overline{\Omega}\circ\mathcal{R}_\gamma\circ \Omega_0^{-1}\circ P(x), \quad\forall\gamma\in \text{PSL}(2,\mathbb{R}).
\end{equation*}
This definition coincides with the homogeneous action of $\text{PSL}(2,\mathbb{R})$ on $F(M)$. Following the discussion of Lemma 3.2 of \cite{LN2021}, we conclude the following result. 

\begin{Prop}\label{prop_conv}
Given any sequence of laminar measures $\mu_i$ on $\mathcal{F}(M,\frac{1}{i})$, the sequence of induced measures $\overline{\Omega}_*\mu_i$ on $F(M)$ converges weakly to a probability measure $\nu$, then $\nu$ is invariant under the homogeneous action of $\text{PSL}(2,\mathbb{R})$.
\end{Prop}

Let $G<\text{PSL}(2,\mathbb{R})$ be a Fuchsian subgroup, then $\mathbb{H}^2/G$ is a closed hyperbolic surface with genus $\geq 2$ whose fundamental domain is represented by $U$.
And let $\phi\in\mathcal{F}(M,\epsilon)$ equivariant with respect to a representation $\Pi$ of $G$ in $\pi_1(M)<SO(n,1)$. The image of $\phi(\mathbb{H}^2)$ in $M$ is closed minimal surface in $M$ whose fundamental group is $\Pi$. We define a laminar measure associated with $\phi$ as follows.  \begin{equation}\label{def_delta_phi}
    \delta_\phi(f)=\frac{1}{vol(U)}\int_U f(\phi\circ\gamma)d\nu_0(\gamma), \quad\forall f\in C^0(\mathcal{F}(M,\epsilon)),
\end{equation}
where $\nu_0$ denotes the bi-invariant measure on $\text{PSL}(2,\mathbb{R})$. 

\subsection{Equidistribution}
In this section, we assume the dimension of $M$ is odd.
Adapting the methods of Proposition 6.1 of \cite{LN2021} and Theorem 5.7 in \cite{Labourie}, we prove the following result. 

\begin{Prop}\label{prop_measure_convergence}
For any $i\in\mathbb{N}$, there is a lamination $\phi_i$ in $\mathcal{F}(M,\frac{1}{i})$ equivariant with respect to a representation of a Fuchsian group $G_i<\text{PSL}(2,\mathbb{R})$ in $\pi_1(M)$, such that $\overline{\Omega}_*\delta_{\phi_i}$ converges to the Lebesgue measure $\mu_{Leb}$ on $F(M)$ as $i\rightarrow\infty$. 
\end{Prop}

\begin{proof}[Sketch of the proof]
Let $\Tilde{\mathcal{T}}$ be the space of tripods $\Tilde{X}=(x_1,x_2,x_3)$, where $x_1,x_2,x_3\in\partial_\infty\mathbb{H}^n$. Each element $\Tilde{X}$ determines an ideal triangle $\Delta(\Tilde{X})$ in $\mathbb{H}^n$. Let $b(\Tilde{X})$ be the barycenter of $\Delta(\Tilde{X})$. Denote by $(e_1(\Tilde{X}),e_2(\Tilde{X}))$ the orthonormal basis of $\Delta(\Tilde{X})$, and denote by $(e_3(\Tilde{X},\cdots,e_n(\Tilde{X}))$ the orthonormal basis of the normal bundle of $\Delta(\Tilde{X})$ in $\mathbb{H}^n$. Thus, each tripod $\Tilde{X}$ determines a point $(b(\Tilde{X}), e_1(\Tilde{X}),\cdots,e_n(\Tilde{X}))$ in $F(\mathbb{H}^n)$, which represents the frame bundle of $\mathbb{H}^n$. 

Consider the closed manifold $M=\mathbb{H}^n/\pi_1(M)$. Let $X$ be the corresponding point of $\Tilde{X}$ in $\mathcal{T}:=\Tilde{\mathcal{T}}/\pi_1(M)$, and let $F(X)$ be the point in the frame bundle $F(M)$ corresponding to $(b(\Tilde{X}), e_1(\Tilde{X}),\cdots,e_n(\Tilde{X}))$ in $F(\mathbb{H}^n)$.

A quintuple $(X,Y,l_1,l_2,l_3)$ is called a \emph{triconnected pair of tripods} (Definition 10.1.1 of \cite{Kahn-Labourie-Mozes}) if $X,Y\in\mathcal{T}$ and $l_1,l_2,l_3$ are three distinct homotopy classes of paths connecting $X$ to $Y$. The space of triconnected pair of tripods is denoted by $\mathcal{TT}$.
Let $\pi^1$ and $\pi^2$ be the forgetting maps from $\mathcal{TT}$ to $\mathcal{T}$\begin{equation*}
    \pi^1:(X,Y,l_1,l_2,l_3)\mapsto X, \quad \pi^2:(X,Y,l_1,l_2,l_3)\mapsto Y.
\end{equation*}
Moreover, let $\overline{\pi^1}$ and $\overline{\pi^2}$ be the corresponding maps from $\mathcal{TT}$ to $F(M)$ \begin{equation*}
    \overline{\pi^1}:(X,Y,l_1,l_2,l_3)\mapsto F(X), \quad \overline{\pi^2}:(X,Y,l_1,l_2,l_3)\mapsto F(Y).
\end{equation*}

In addition, there is a weighted measure $\mu_{\epsilon,R}$ on $\mathcal{TT}$ as defined in Definition 11.2.3 of \cite{Kahn-Labourie-Mozes}. If $(X,Y,l_1,l_2,l_3)$ is in the support of $\mu_{\epsilon,R}$, then the ideal triangles determined by $X$ and $Y$ can be glued to an \emph{$(\epsilon,R)$-almost closing pair of pants} (see Definition 9.1.1 of \cite{Kahn-Labourie-Mozes}). Moreover, it follows from the mixing property that for fixed $\epsilon$, as $R\rightarrow\infty$, $\overline{\pi^1}_*\mu_{\epsilon,R}$ and $\overline{\pi^2}_*\mu_{\epsilon,R}$ both converge to the Lebesgue measure $\mu_{Leb}$ on $F(M)$.

Arguing like Theorem 5.7 of \cite{Labourie}, we can choose a sequence $R_j\rightarrow\infty$ as $j\rightarrow\infty$, and a sequence of measures $\mu_{\frac{1}{j},R_j}$. Then we approximate each $\mu_{\frac{1}{j},R_j}$ by another weighted measure $\nu_j$ supported in finitely many pleated pair of pants $P^1_j,\cdots,P^{N_j}_j$, which can be glued together to get essential surfaces $\Sigma_j^1,\cdots,\Sigma_j^{M_j}$ in $M$.
When $j$ is sufficiently large and for each $1\leq k\leq M_j$, $\Sigma_j^k$ is $(1+\frac{1}{j})$-quasigeodesic,
and the projection from $\Sigma_j^k$ to the unique minimal surface $S_j^k$ homotopic to $\Sigma_j^k$ is $(1+\frac{1}{j})$-bi-Lipschitz and it has distance uniformly bounded by $O(\frac{1}{j})$. For this reason, we can further approximate $\nu_j$ by a weighted measure supported in $S_j^1\cup\cdots\cup S_j^{M_j}$.
$S_j^k$ is obtained by a lamination $\phi_j^k\in\mathcal{F}(M,\frac{1}{j})$, in fact, it's the image of $\phi_j^k(\mathbb{H}^2)$ in $M$, and thus associated with the laminar measure $\delta_{\phi_j^k}$, we have the following lemma.

\begin{Lemma}
For any $j\in\mathbb{N}$, there exist a finite sequence of laminations $\phi_j^1,\cdots,\phi_j^{M_j}$ in $\mathcal{F}(M,\frac{1}{j})$, and $\theta_j^1,\cdots,\theta_j^{M_j}\in (0,1)$ with $\theta_j^1+\cdots +\theta_j^{M_j}=1$, such that each $\phi_j^k$ is equivariant with respect to a representation of a Fuchsian group in $\pi_1(M)$, and the laminar measure \begin{equation*}
    \mu_j=\sum_{k=1}^{M_j}\theta_j^k\delta_{\phi_j^k}
\end{equation*} 
satisfies that $\overline{\Omega}_*\mu_j$ converges to the Lebesgue measure $\mu_{Leb}$ on $F(M)$ as $j\rightarrow\infty$. 
\end{Lemma}

Next,  for $2\leq l\leq n-1$, we define \begin{align*}
    \mathcal{P}_l:=\{F(P)\subset F(M), &\text{ where }P\text{ is a $l$-dimensional closed totally geodesic}\\
    &\text{ submanifold of }M\}.
\end{align*}
Then $\mathcal{P}:\underset{l=2}{\overset{n-1}{\cup}}\mathcal{P}_l$ contains at most countably many candidates.
Therefore, we can find a decreasing sequence of tubular neighborhoods $\{B_k\}\subset F(M)$, so that for any $k\in\mathbb{N}$, $B_k$ covers $\mathcal{P}$ and it satisfies  $\mu_{Leb}(B^k)<2^{-2k-1}$ and $\mu_{Leb}(\partial B_k)=0$. 
In consequence of previous lemma, after passing to a subsequence, we have $\overline{\Omega}_*\mu_j(B_k)<2^{-2k}$. Additionally, as argued in Lemma 6.2 of \cite{LN2021}, we can find a subsequence $\{j_i\}$, and $\phi_i\in\{\phi_{j_i}^1,\cdots,\phi_{j_i}^{M_{j_i}}\}$, 
such that $\overline{\Omega}_*(\delta_{\phi_i})(B_k)<2^{-k}$.

As a result of Proposition \ref{prop_conv}, as $i\rightarrow\infty$, $\overline{\Omega}_*\delta_{\phi_i}$ converges weakly to a probability measure $\nu$ on $F(M)$. $\nu$ is invariant under the homogeneous action of $\text{PSL}(2,\mathbb{R})$, and it satisfies that \begin{equation}\label{nu_ineq}
    \nu(B_k)<2^{-k}.
\end{equation}
To finish the proof, we need the following lemma.

\begin{Lemma}
$\nu=\mu_{Leb}$.
\end{Lemma}

\begin{proof}
According to the ergodic decomposition theorem (\cite{Hasselblatt-Katok}), $\nu$ can be expressed by a linear combination of the ergodic measures for $\text{PSL}(2,\mathbb{R})$-action on $F(M)$. Moreover, Ratner's measure classification theorem (see \cite{Ratner} or \cite{Shah}) says that any ergodic $\text{PSL}(2,\mathbb{R})$-invariant measure on $F(M)$ is either an invariant probability measure supported on a finite union of $\{P_k\}\subset\mathcal{P}$, or it is identical to $\mu_{Leb}$. Thus, we can write $\nu$ as \begin{equation*}
    \nu=a_1\mu_{Leb}+a_2\mu_{\mathcal{P}_2}+\cdots +a_{n-1}\mu_{\mathcal{P}_{n-1}},
\end{equation*}
where $a_1+a_2+\cdots +a_{n-1}=1$ and $\mu_{\mathcal{P}_l}$ represents an ergodic measure supported on $\mathcal{P}_l$, $2\leq l\leq n-1$. 
By (\ref{nu_ineq}), for all $k\in\mathbb{N}$, \begin{equation*}
    a_2+\cdots+a_{n-1}=a_2\mu_{\mathcal{P}_2}(B_k)+\cdots +a_{n-1}\mu_{\mathcal{P}_{n-1}}(B_k)\leq \nu(B_k)<2^{-k}.
\end{equation*}
So \begin{equation*}
    a_1=1-a_2-\cdots-a_{n-1}>1-2^{-k},\quad\forall k\in\mathbb{N}.
\end{equation*}
We must have $a_1=1$, and therefore $\nu=\mu_{Leb}$.
\end{proof}
Proposition \ref{prop_measure_convergence} follows immediately from the lemma.
\end{proof}

\subsection{Average area ratio formula}
\begin{Lemma}[average area ratio formula]\label{lemma_area ratio formula}
Let $(N,g)$ be a closed Riemannian manifold that also has odd dimension $n$, and let $F$ be a smooth map that takes $(N,g)$ to $(M,h_0)$. 
For $i\in\mathbb{N}$, we pick a lamination $\phi_i\in\mathcal{F}(M,\frac{1}{i})$ equivariant with respect to a representation of $G_i<\text{PSL}(2,\mathbb{R})$ in $\pi_1(M)$, and it satisfies Proposition \ref{prop_measure_convergence}. Let $S_i$ be the image of $\phi_i(\mathbb{H}^2)$ in $M$.
Then we have \begin{equation*}
    \text{Area}_{F}(g/h_0)=\lim_{i\rightarrow\infty}\frac{\text{area}_g(F^{-1}(S_i))}{4\pi(g_i-1)}.
\end{equation*}
\end{Lemma}

\begin{proof}
Recall that $|\Lambda^2F|^{-1}_g$ is a function defined almost everywhere on $Gr_2(M)$.
Since $|\Lambda^2F|_g$ can be regarded as a smooth function on $F(M)$ by \begin{equation*}
    |\Lambda^2F|_g:F(M)\rightarrow\mathbb{R},\quad (x,e_1,e_2\cdots,e_n)\mapsto |\Lambda^2F|_g(x,\text{span}(e_1,e_2)),
\end{equation*}
based on the definition (\ref{def_jac}), $|\Lambda^2F|^{-1}_g$ is also seen as a function defined almost everywhere on $F(M)$. Thus, Proposition \ref{prop_measure_convergence} implies that
\begin{equation*}
    \text{Area}_F(g/h_0)=\mu_{Leb}(|\Lambda^2F|^{-1}_g)=\lim_{i\rightarrow\infty}\overline{\Omega}_*\delta_{\phi_i}(|\Lambda^2F|^{-1}_g).
\end{equation*}
In light of the definition of laminar measure $\delta_{\phi_i}$ in (\ref{def_delta_phi}), we have \begin{equation*}
    \overline{\Omega}_*\delta_{\phi_i}(|\Lambda^2F|^{-1}_g) =\frac{1}{vol(U_i)}\int_{U_i} |\Lambda^2F|^{-1}_g\circ\overline{\Omega}(\phi_i\circ\gamma)d\nu_0(\gamma),
\end{equation*}
where $U_i$ is the fundamental domain of $\text{PSL}(2,\mathbb{R})/G_i$. Set $x=\gamma(i)$. Since the hyperbolic surface $\mathbb{H}^2/G_i$ has area equal to $4\pi(g_i-1)$, where $g_i\geq 2$ denotes the genus. The above expression also can be written as 
\begin{align*}
    \overline{\Omega}_*\delta_{\phi_i}(|\Lambda^2F|^{-1}_g) &= \frac{1}{4\pi(g_i-1)}\int_{\mathbb{H}^2/G_i} |\Lambda^2F|^{-1}_g(\phi_i(x),(d\phi_i)_xT_x\mathbb{H}^2)dA_{h_{\mathbb{H}^2/G_i}}(x)\\
    &= \frac{1}{4\pi(g_i-1)}\int_{S_i}\frac{|\Lambda^2F|^{-1}_g(y,T_yS_i)}{C_i^2(\phi_i^{-1}(y))}dA_{h_0}(y),
\end{align*}
where $C_i^2>0$ denotes the conformal factor between the hyperbolic metric on $\mathbb{H}^2 /G_i$ and the pull-back metric of $h_0$ by $\phi_i$, namely $\phi_i^*(h_0)=C_i^2h_{\mathbb{H}^2/G_i}$. Since the Gaussian curvature on $S_i$ has the form $-1-\frac{1}{2}|A|^2(x)$, we have \begin{equation*}
     1\leq \frac{1}{C_i^2}\leq 1+\frac{1}{2}|A|^2_{L^\infty(S_i)}.
\end{equation*}
On the other hand, the co-area formula yields that \begin{equation*}
    \int_{S_i} |\Lambda^2F|^{-1}_g(y,T_yS_i)dA_{h_0}(y)=\text{area}_g(F^{-1}(S_i)).
\end{equation*}
Combining these formulas, we have \begin{equation*}
    \frac{\text{area}_g(F^{-1}(S_i))}{4\pi(g_i-1)}\leq \overline{\Omega}_*\delta_{\phi_i}(|\Lambda^2F|^{-1}_g)\leq (1+\frac{1}{2}|A|^2_{L^\infty(S_i)})\frac{\text{area}_g(F^{-1}(S_i))}{4\pi(g_i-1)}.
\end{equation*}
Since $|A|^2_{L^\infty(S_i)}\rightarrow 0$ as $i\rightarrow\infty$ (see \cite{Jiang}), the lemma follows immediately from the squeeze theorem.
\end{proof}

\section{Gromov's conjecture in dimension three}\label{section_gromov_3}
In this section, we discuss the proof of Theorem \ref{Thm_gromov_3}.
First of all, the fact that $\text{deg}\,F=1$ follows from Corollary 0.3 of \cite{Wang1993} and the geometrization theorem for 3-manifolds. Next, it's easy to see, the induced map $F_*:\pi_1(N)\rightarrow \pi_1(M)$ is surjective, since otherwise, it factors through a $d$-fold covering space of $M$ with $d>1$, and thus $\text{deg}\,F\geq d>1$, violating the degree one observation. In addition, $\pi_1(M)$ is a Hopfian group (for example, see 15.13 of \cite{Hempel}), so the surjectivity of $F_*$ can be upgraded to be an isomorphism, which makes $F$ a homotopy equivalence between $N$ and $M$ due to Whitehead's theorem. Furthermore, the Mostow rigidity theorem indicates that $F$ is homotopic to an isometry. For this reason, we can simplify the conditions of Theorem \ref{Thm_gromov_3} as follows.

\subsection{A simpler version of Theorem \ref{Thm_gromov_3}}
\begin{namedtheorem}[Theorem 1.4']\label{theorem2'}
Let $(M,h_0)$ be a closed hyperbolic 3-manifold. There exists a small neighborhood $\mathcal{U}$ of $h_0$ in the $C^2$-topology, such that for all Riemannian metric $h\in\mathcal{U}$ with $R_h\geq -6$, and for any smooth map $F:(M,h)\rightarrow (M,h_0)$, it has $\text{deg}\,F=1$ and it is homotopic to the identity, 
we have \begin{equation*}
   \quad \text{Area}_F(h/h_0)\geq 1.
\end{equation*}
Moreover, the equality holds if and only if $F$ is an isometry between $h$ and $h_0$.
\end{namedtheorem}

Let $S_i$ be the minimal surface in $M$ with respect to $h_0$ defined in Lemma \ref{lemma_area ratio formula}. 
The inverse $F^{-1}(S_i)$ is also a closed surface in $M$, but note that $F^{-1}(S_i)$ is not necessarily homotopic to $S_i$. In fact, we can only find the following relation of their genus. The Gromov norms of $S_i$ and $F^{-1}(S_i)$ satisfies that \begin{equation*}
    |\deg F|||S_i||\leq ||F^{-1}(S_i)||.
\end{equation*}
Here $\deg F=1$. And for any closed surface $S$ with genus $g(S)$, $||S||=\frac{4\pi(g(S)-1)}{v_2}$, where $v_2$ is a fixed number representing the supreme area of geodesic 2-simplices in $\mathbb{H}^2$. As an immediate result, we have $g(F^{-1}(S_i))\geq g(S_i)$.

To compare the areas of surfaces with respect to the induced metric of $h$ in different homotopic classes, we hope to find a global area-minimizing surface. In general, the existence and the topology of such a surface are complicated. But if there is a minimal surface with suitable curvature conditions, then adapting Uhlenbeck's method in \cite{Uhlenbeck}, we can check the uniqueness of a closed minimal surface of any type, which is the key point of the proof.

\subsection{Proof of \nameref{theorem2'}}
For each $i\in\mathbb{N}$, let $\Tilde{M}_i$ be the covering space of $M$ such that $\pi_1(\Tilde{M}_i)\cong \pi_1(S_i)$.
Let $\Tilde{F}_i$ be the corresponding lift of $F$ that maps $\Tilde{F}_i^{-1}(\Tilde{M}_i)\simeq\Tilde{M}_i$ to $\Tilde{M}_i$.
The lift of $S_i$ in $\Tilde{M}_i$ still has fundamental group $\pi_1(S_i)$, so we denote it by $S_i$ as well.
By assumption, $F$ is homotopic to identity, thus there is a continuous map $H:M\times [0,1]\rightarrow M$ with $H(x,0)=x$ and $H(x,1)=F(x)$ for any $x\in M$. Since $M$ is compact, the length of the path of $H$ between $x$ and $F(x)$ is uniformly bounded by a constant $C>0$. Now let $\Tilde{H}_i$ be the lift of $H$ that connects $\Tilde{F}_i$ to the identity map on $\Tilde{M}_i$. For all $y\in\Tilde{M}_i$, the length of the path between $y$ and $\Tilde{F}_i(y)$ is therefore uniformly bounded by the same constant $C$. 
So $\Tilde{F}_i$ is proper, meaning $\Tilde{F}^{-1}_i(S_i)$ is a closed set, and therefore it is a $k$-fold cover of $F^{-1}(S_i)$ for some finite number $k$. If $k>1$, the image of $\Tilde{F}^{-1}_i(S_i)$ under $\Tilde{F}_i$ is either a closed surface with Euler characteristic equal to $k\chi(S_i)$, and therefore having genus $kg(S_i)-k+1>g(S_i)$, or it is a union of at least two surfaces with genus $\geq g(S_i)$. However, both cases are impossible because the image cannot be identified with $S_i\subset\Tilde{M}_i$. 
We must have $k=1$.
Consequently, 
the covering map from $\Tilde{F}_i^{-1}(S_i)$ to $F^{-1}(S_i)$ is one-to-one. 

On the other hand, the classical result \cite{SY} verifies the existence of area-minimizing surface $\Sigma_i\subset (M,h)$ in the homotopy class of $S_i$.
And based on Theorem 4.3 of \cite{Lowe}, there exist a $C^2$-neighborhood $\mathcal{U}_0$ of $h_0$ and $N_0\in\mathbb{N}$, so that when $h\in\mathcal{U}_0$ and $i\geq N_0$, $\Sigma_i$ is the unique minimal surface in $(M,h)$ homotopic to $S_i$.
Furthermore, let $D_i$ ($\Omega_i$) be the lifts of $S_i$ ($\Sigma_i$, respectively) in $B^3$. These discs $D_i$ and $\Omega_i$ are asymptotic and at a uniformly bounded Hausdorff distance to each other, as $h\rightarrow h_0$, $\Omega_i$ converges uniformly on compact sets to $D_i$ in $C^{2,\alpha}$. Therefore, replacing $\mathcal{U}_0$ by a smaller subset or replacing $N_0$ by a larger integer if necessary, we can assume that if $h\in\mathcal{U}_0$ and $i\geq N_0$, then there exists a smooth map $f_i$ on $D_i$ with $|f_i|_{C^{2,\alpha}}<1$, such that $\Omega_i$ can be represented by a graph of $f_i$ over $D_i$. More precisely, Let $n_i$ be the unit normal vector field of $D_i$, then we have the following diffeomorphism.
\begin{equation*}
    F_i:D_i\rightarrow\Omega_i,\quad F_i(x)=\cosh (f_i(x))x+\sinh (f_i(x))n_i(x).
\end{equation*}

Notice that the minimal disc $\Omega_i$ has mean curvature equal to zero with respect to $h$, so the mean curvature $H_{h_0}(\Omega_i)$ with respect to $h_0$ has a uniform bound determined by the perturbation of $h$ and $\nabla h$. Since $h$ is $C^2$-close to $h_0$, we have \begin{equation}\label{equ_1_gromov3}
    |H_{h_0}(\Omega_i)|_{C^{0,\alpha}}=O(|h-h_0|_{C^2}),\quad\forall i\geq N_0.
\end{equation}
According to the Schauder estimates for elliptic PDE, there exists a constant $c_0>0$, such that for any $i\geq N_0$, \begin{equation}\label{equ_2_gromov3}
    |f_i|_{C^{2,\alpha}}\leq c_0\big(|f_i|_{L^\infty}+|H_{h_0}(\Omega_i)|_{C^{0,\alpha}}\big).
\end{equation} 

Besides, for $i\geq N_0$, suppose the principle curvature of $D_i$ with respect to $h_0$ satisfies that $\underset{i\geq N_0}{\sup}|\lambda(D_i)|_{L^\infty}<1$. Uhlenbeck (\cite{Uhlenbeck}) shows that $\mathbb{H}^3$ is foliated by a sequence of equidistant discs relative to $D_i$. We denote by $D_i^r$ the disc with a fix distance $r$ to $D_i$, it has mean curvature \begin{equation*}
    H_{h_0}(D_i^r)=\frac{2(1-\lambda(D_i)^2)\tanh r}{1-\lambda(D_i)^2\tanh^2 r}.
\end{equation*}
Let $R_i^+$ and $R_i^-$ be the supremum and infimum of $r$ such that $\Omega_i$ meets $D_i^r$, and the intersections points are $x_i^+$, $x_i^-$, respectively. Since $R_i^+,R_i^-\rightarrow 0$ as $h\rightarrow h_0$, we may assume \begin{equation*}
    \min\{\frac{d}{dr}(\tanh r)|_{r=R_i^+}, \,\frac{d}{dr}(\tanh r)|_{r=R_i^-}\}=\min\{\frac{1}{\cosh^2 R_i^+},\,\frac{1}{\cosh^2 R_i^-}\}\geq \frac{1}{2},
\end{equation*}
then we have \begin{align*}
    |H_{h_0}(\Omega_i)|_{L^\infty} &\geq \max\{|H_{h_0}(D_i^{R_i^+}(x_i^+)|,\,|H_{h_0}(D_i^{R_i^-}(x_i^-)|\}\\
    &\geq 2(1-|\lambda(D_i)|^2_{L^\infty})\max\{|\tanh R_i^+|,\,|\tanh R_i^-|\}\\
    &\geq (1-|\lambda(D_i)|^2_{L^\infty})\max\{|R_i^+|,\,|R_i^-|\}.
\end{align*}
Since $\Omega_i$, described by the graph $f_i$, is bounded between $D_i^{R_i^+}$ and $D_i^{R_i^-}$, the above result indicates the existence of a uniform constant $c_1>0$, such that \begin{equation}\label{equ_3_gromov3}
    |f_i|_{L^\infty}\leq c_1|H_{h_0}(\Omega_i)|_{L^\infty},\quad \forall i\geq N_0.
\end{equation}
Combining (\ref{equ_1_gromov3})-(\ref{equ_3_gromov3}), we obtain \begin{equation*}
    |f_i|_{C^{2,\alpha}}=O(|h-h_0|_{C^2}),\quad \forall i\geq N_0.
\end{equation*}
And therefore, the principal curvatures of $\Sigma_i$ with respect to $h_0$ and $h$ satisfy that \begin{align*}
   & |\lambda_{h_0}(\Sigma_i)|^2_{L^\infty}=O(|h-h_0|^2_{C^2}),\quad \forall i\geq N_0,\\
   \Longrightarrow &\, |\lambda_{h}(\Sigma_i)|^2_{L^\infty}=|\lambda_{h_0}(\Sigma_i)|^2_{L^\infty}+O(|h-h_0|^2_{C^2})=O(|h-h_0|^2_{C^2}),\quad \forall i\geq N_0.
\end{align*}
Clearly, the sectional curvature of $(M,h)$ has the property \begin{equation*}
    |K_h|_{L^\infty}=-1+O(|h-h_0|_{C^2}),\quad \forall i\geq N_0.
\end{equation*}
Thus, we can find $\mathcal{U}\subset\mathcal{U}_0$ and $N\geq N_0$, such that if $h\in\mathcal{U}$ and $i\geq N$, then the principal curvatures of $\Sigma_i$ with respect to $h$ and the sectional curvature of $(M,h)$ satisfy that
\begin{equation*}
    |\lambda_h(\Sigma_i)|_{L^\infty}^2<-\sup K_h.
\end{equation*}
In the lemma below, we apply Uhlenbeck's method \cite{Uhlenbeck}, as well as the comparison result associated with Riccati equations, to prove that $\Sigma_i$ is the unique closed minimal surface in $(\Tilde{M}_i,h)$, thus minimizing the area among all closed surfaces.

\begin{Lemma}
Let $\Sigma$ be a minimal surface in $(M,h)$ whose fundamental group injectively includes in $\pi_1(M)$, the principal curvatures of $\Sigma$ with respect to $h$, denoted by $\pm\lambda$, and the sectional curvature $K$ of $(M,h)$ satisfy that \begin{equation}\label{lambda_assumption}
    |\lambda(\Sigma)|^2_{L^\infty}<-\sup K.
\end{equation}
Let $\Tilde{M}$ be the cover of $M$ with $\pi_1(\Tilde{M})\cong \pi_1(\Sigma)$. Then $\Sigma$ is the unique closed minimal surface in $\Tilde{M}$.
\end{Lemma}

\begin{proof}
Denote the supreme of the sectional curvature on $(M,h)$ by $-k^2$, where $k>0$. 
Let $f(x)$ be the distance function from a fixed point in $\Tilde{M}\setminus\Sigma$ to $x\in \Sigma$.
By (\ref{lambda_assumption}), for any $X\in T_x\Sigma$, \begin{equation}
    D^2f(X,X)=\text{Hess}\,f(X,X)- A(X,X)(f)\geq \frac{k}{\tanh(kf(x))}-k>0.
\end{equation}
It turns out that $f$ is a convex function, thus, there's only one critical point that attains the minimum. 
As a result, $\exp|_{N\Sigma}$ maps injectively from the normal bundle $N\Sigma$ to $\Tilde{M}$.

Furthermore, we show that $\exp|_{N\Sigma}$ is a diffeomorphism, and thus $\Tilde{M}$ is foliated by a family of surfaces $\{\Sigma_r\}_{r\in\mathbb{R}}$, where $\Sigma_r$ is the surface at the fixed distance $r$ to $\Sigma$. 
To see this, we introduce some notations beforehand.
For $x\in \Sigma$, choose an oriented, orthonormal basis $\{e_1,e_2\}$ for $T_{x}\Sigma$, and a unit vector $e_3$ for $N_{x}\Sigma$. Then we obtain an orthonormal frame by applying parallel transport along $\exp|_{N\Sigma}$. Since $\Sigma$ is a minimal surface, the principal curvatures satisfy that $\lambda_{1}=-\lambda_{2}:=\lambda$, 
we assume $\lambda\geq 0$ in the following computation. 
Let $V_{i}(r)=v^i(r)e_i(r)$ be the Jacobi field along $\exp(re_3)$, where $i=1,2$, it satisfies that $v^i(0)=1$, $(v^i)'(0)=\lambda_{i}=\pm \lambda$.

On the other hand, let $\overline{\Sigma}$ be a minimal surface in $\Tilde{M}$ with respect to an ambient metric $\overline{h}$ of constant sectional curvature $-k^2$, and its principal curvature satisfies that $\overline{\lambda}=\lambda$. We do not require the existence of $\overline{\Sigma}$, it's only used for comparison in the computation. Similar to the notations defined above, let $\overline{e}_1,\cdots,\overline{e}_3$ be the corresponding frame on $\Tilde{M}$ with respect to $\overline{h}$, and let $\overline{V}_{i}(r)=\overline{v}^i(r)\overline{e}_i(r)$ be the Jacobi field along $\exp(r\overline{e}_3)$ which shares the same initial data with $V_{i}(r)$. 
Since $\lambda<k$, we have \begin{equation*}
    \overline{v}^i(r)=
    \cosh{kr}\pm \dfrac{\lambda\sinh{kr}}{k}>0,\quad i=1,2.
\end{equation*}
From \begin{equation*}
    (v^{i})''=-R(e_3,e_i,e_i,e_3)v^{i}\geq k^2\, \overline{v}^{i}=(\overline{v}^{i})''
\end{equation*}
and the initial data,
the graph of $v^i$ lies above that of $\overline{v}^i$, thus above the horizontal axis. The non-vanishing Jacobi fields ensure that the induced metric on $\Sigma_r$ in $(\Tilde{M},h)$ is nonsingular for all $r\in\mathbb{R}$. In addition, we've seen that $\exp|_{N\Sigma}: N\Sigma\rightarrow \Tilde{M}$ is injective, and therefore also bijective, so it is a diffeomorphism and $\Tilde{M}$ admits a foliation structure.

Next, let $\lambda_{i}(r)\,(i=1,2)$ be the principal curvatures on $\Sigma_r$, and denote by $\overline{\lambda}_{i}(r)\,(i=1,2)$ the principal curvatures of the $r$-equidistant surface to $\overline{\Sigma}$ with respect to $\overline{h}$. Notice that each $\lambda_{i}(r)$ satisfies the Riccati equation \begin{equation*}
    \lambda_{i}'(r)=\lambda^2_{i}(r)+R(e_3,e_i,e_i,e_3)(r).
\end{equation*}
Then it follows from the comparison theorem associated with Riccati equations (for instance, see Theorem 3.1 of \cite{Warner}) that \begin{align*}
    \lambda_{1}(r) &\geq\overline{\lambda}_{1}(r)= k\,\frac{k\tanh(kr)+\lambda}{k+\lambda\tanh(kr)}>0,\\
    \lambda_{2}(r) &\geq \overline{\lambda}_{2}(r)= k\,\frac{k\tanh(kr)-\lambda}{k-\lambda\tanh(kr)}.
\end{align*}
It follows from $\lambda^2\leq k^2$ that \begin{equation*}
    \lambda_{1}(r)+\lambda_{2}(r) \geq 2k\,\frac{(k^2-\lambda^2)\tanh(kr)}{k^2-\lambda^2\tanh^2(kr)}>0.
\end{equation*}
Therefore, for any $r\in\mathbb{R}$, $\Sigma_r$ is strictly mean convex with respect to the metric induced by $h$.

Finally, we prove the uniqueness. Assume that $\Sigma'$ is another closed minimal surface in $(\Tilde{M},h)$, and let $R_+$ and $R_-$ be the supremum and infimum of $r$ such that $\Sigma'$ intersects $\Sigma_r$, respectively, then $R_+$ and $R_-$ are both finite. However, due to the maximum principle, $\Sigma'$ cannot be tangential to any strictly mean convex slice $\Sigma_r$ with $r\neq 0$. Therefore, we must have $\Sigma'\subset\Sigma_0=\Sigma$.
\end{proof}

Now we finish the proof of Theorem \ref{Thm_gromov_3}. From the previous lemma, when $i\in\mathbb{N}$ is sufficiently large, $\Sigma_i$ is the area-minimizer among all closed surfaces in $\Tilde{M}_i$ with respect to the induced metric of $h$, it yields that \begin{equation*}
    \text{area}_h(F^{-1}(S_i))=\text{area}_h(\Tilde{F}_i^{-1}(S_i))\geq \text{area}_h(\Sigma_i).
\end{equation*} 
Combining it with the area comparison in Theorem 5.1 of \cite{LN2021}, we have \begin{equation*}
    \text{Area}_F(h/h_0)=\lim_{i\rightarrow\infty}\frac{\text{area}_h(F^{-1}(S_i))}{4\pi(g_i-1)}\geq\lim_{i\rightarrow\infty} \frac{\text{area}_h(\Sigma_i)}{4\pi(g_i-1)}\geq 1.
\end{equation*}


Moreover, when the equality holds, it follows from the equality of Theorem 5.1 of \cite{LN2021} that $h=F^*(h_0)$,  $F$ is an isometry between $h$ and $h_0$.

\section{Proof of Theorem \ref{Thm_odd}}\label{section_odd}
In the end, we prove Theorem \ref{Thm_odd} in this last section.
\subsection{Proof of inequality}
The proof follows directly from \cite{LN2021}, but for readers' convenience, it is stated as follows.

We let \begin{equation*}
    \alpha:
    =\lim_{i\rightarrow\infty}\frac{\text{area}_h(\Pi_i)}{4\pi(g_i-1)}.
\end{equation*}
For any $\delta>0$, we can take $i$ sufficiently large,  such that \begin{equation*}
    \lim_{i\rightarrow\infty}\frac{\text{area}_h(\Pi_i)}{4\pi(g_i-1)}<\alpha+\delta.
\end{equation*}
Let $\Pi_i^k$ be the $k$-cover of $\Pi_i$. Since $\Pi_i$ has genus $g_i\geq 2$, \begin{equation*}
    4\pi(g_i^k-1)\geq k\,4\pi(g_i-1),
\end{equation*}
then the least area surface in the homotopy class of $\Pi_i^k$ with respect to $h$ satisfies that \begin{equation}\label{area_of_pi_k}
    \lim_{i\rightarrow\infty}\frac{\text{area}_h(\Pi_i^k)}{4\pi(g_i^k-1)}\leq \lim_{i\rightarrow\infty}\frac{k\,\text{area}_h(\Pi_i)}{k\,4\pi(g_i-1)}
    =\lim_{i\rightarrow\infty}\frac{\text{area}_h(\Pi_i)}{4\pi(g_i-1)}
    <\alpha+\delta.
\end{equation}
According to M{\"u}ller-Puchta's formula (see \cite{Kahn-Markovic}), there exists a constant $c_i$ that depends only on $M$ and $i$, such that the following is true when $g_i$ is large. \begin{equation}\label{Muller-Puchta}
    s(M,g_i^k,\frac{1}{i})\geq (c_ig_i^k)^{2g_i^k}.
\end{equation} 
Define $L_i^k$ in the following way \begin{equation}\label{L_i^k}
    4\pi(L_i^k-1)=(\alpha+\delta)4\pi(g_i^k-1)\Longrightarrow \lim_{k\rightarrow\infty}\frac{g_i^k}{L_i^k}=\frac{1}{\alpha+\delta}.
\end{equation}
Combining (\ref{Muller-Puchta}) with (\ref{area_of_pi_k}), we have that \begin{equation*}
    \#\{\text{area}_h(\Pi)\leq 4\pi(L_i^k-1):\Pi\in S_{\frac{1}{i}}(M)\}\geq (c_ig_i^k)^{2g_i^k},
\end{equation*}
Therefore, by (\ref{L_i^k}),\begin{align*}
     E(h)&=\underset{i\rightarrow \infty}{\lim}\,\underset{k\rightarrow\infty}{\lim\inf}\,\dfrac{\ln\#\{\text{area}_h(\Pi)\leq 4\pi(L_i^k-1):\Pi\in S_{\frac{1}{i}}(M)\}}{L_i^k\ln L_i^k}\\
     &\geq \underset{i\rightarrow \infty}{\lim}\,\underset{k\rightarrow\infty}{\lim\inf}\,\dfrac{(c_ig_i^k)^{2g_i^k}}{L_i^k\ln L_i^k}\\
     &=\frac{2}{\alpha+\delta}.
\end{align*}
Since $\delta$ is an arbitrarily small positive number, we conclude that \begin{equation}\label{area_e(h)}
     \text{Area}_{\text{Id}}(h/h_0)E(h)=\lim_{i\rightarrow\infty}\frac{\text{area}_h(S_i)}{4\pi(g_i-1)}\,E(h)
     \geq \lim_{i\rightarrow\infty}\frac{\text{area}_h(\Pi_i)}{4\pi(g_i-1)}\,E(h)\geq 2.
\end{equation}

\subsection{Proof of rigidity}
If $\text{Area}_{\text{Id}}(h/h_0)E(h)=2$, then (\ref{area_e(h)}) yields that\begin{equation}\label{limit_1}
    \lim_{i\rightarrow\infty}\frac{\text{area}_h(\Pi_i)}{\text{area}_h(S_i)}=1.
\end{equation}

To make use of this equality, we run the mean curvature flow in $(B^n,h)$ with initial condition $D_i$, which is the lift of $S_i$ in $\mathbb{H}^n$, then we estimate the decay rate of the area. Firstly of all, we need to review and establish some tools for complete, noncompact surfaces moving by mean curvature. The classical short time existence theorem for compact manifolds moving by mean curvature is well-known \cite{Hamilton}.  However, the general theory for complete, noncompact manifolds has not been established in the literature. There are only several essential contributions in some special cases: Ecker-Huisken \cite{Ecker-Huisken} proved the codimension one case in which only a local Lipschitz condition on the initial hypersurface was required. For higher codimensions, Chau-Chen-He \cite{Chau} discussed the case of nonparametric mean curvature flow for flat metrics. The result related to our case is listed as follows.

\begin{Lemma}\label{short_time_existence}
There exist $T>0$ and $C>0$ depending only on $M$, so that for sufficiently large $i\in\mathbb{N}$, we can find a solution $D_i(t)$ to the mean curvature flow in $(M,h)$ with initial condition $D_i(0)=D_i$, where $0\leq t\leq T$. Additionally, the mean curvature of $D_i(t)$ and its derivative are both bounded uniformly by $C$. 
\end{Lemma}
\begin{proof}

Notice that after passing to a subsequence, $D_i$ converges smoothly on compact sets to a disc $D$, and each of them is a cover of a compact surface in $M$. Take $x\in D$, the standard theory indicates that there is a number $T_0>0$, such that for any $k\in\mathbb{N}$, we can find a solution $D_i^k(t)$ to the mean curvature flow with initial condition $\overline{B(x,k)}\cap D_i$, where $0\leq t\leq T_0$. Since $T_0$ depends only on the second fundamental form of $D$, in particular, it's independent of $i$ and $k$. 


Next, in order to apply the Arzela-Ascoli theorem and estimate the mean curvature of $D_i(t)$ and its derivative for any small time $t$, we need the following preparation.
Claim that for any $\delta>0$, and any spacetime $X_i^k=(x_i^k,t)$ of $D_i^k(t)$,  there exists an open neighborhood $U_i^k$ of $X_i^k$, so that the \emph{Guassian density ratio} \begin{equation*}
    \Theta(D_i^k(t),X_i^k,r):=\int_{y\in D_i^k(t-r^2)}\frac{1}{4\pi r^2}\exp\big(-\frac{|y-x_i^k|^2}{4r^2}\big)\,d\mathcal{H}^2(y)
\end{equation*}
satisfies that \begin{equation}\label{density}
    \Theta(D_i^k(t), X_i^k, r)\leq 1+\delta, \quad \forall\, 0<r<d(X_i^k,U_i^k).
\end{equation}
If this wasn't true for some integers $i$ and $k$, then we could pick a sequence $\lambda_j\rightarrow \infty$ as $j\rightarrow\infty$, and \begin{equation*}
    \Theta(\mathcal{D}_{\lambda_j}(D_i^k(t)-X_i^k), 0, \lambda_jr)> 1+\delta,
\end{equation*}
where $\mathcal{D}_\lambda$ denotes the parabolic dilation $\mathcal{D}_{\lambda}(y,t)=(\lambda y,\lambda^2 t)$. 
Since the second fundamental form satisfies $|A|^2_{L^{\infty}(\mathcal{D}_{\lambda_j}(D_i^k(t))-X_i^k)}\rightarrow 0$ as $j\rightarrow\infty$,  $\mathcal{D}_{\lambda_j}(D_i^k(t)-X_i^k)$ converges smoothly to a disc $\bar{D}_i^k$ whose second fundamental form vanishes. However, the inequality above implies that \begin{equation*}
    \lim_{j\rightarrow\infty}\Theta(\bar{D}_i^k,0,\lambda_jr)>1,
\end{equation*} which contradicts the topology of $\bar{D}_i^k$.

We've seen that (\ref{density}) holds, so due to the local regularity theorem in \cite{White}, there is a uniform constant $C_0$ that is independent of $i$ and $k$, so that at any spacetime $X_i^k=(x_i^k,t)$, \begin{equation}\label{white_regularity}
    |A|^2(X_i^k)d(X_i^k,U_i^k)\leq C_0.
\end{equation}
Therefore, Arzela-Ascoli theorem (see page 1494 of \cite{White}) implies the short-time existence of the mean curvature flow with noncompact initial condition $D_i(0)=D_i$ on time interval $[0,T_0]$. 
Moreover, the interior estimate (\ref{white_regularity}) validates the condition of the maximum principle (\cite{Ecker-Huisken}, Theorem 4.3). Arguing like Theorem 4.4 of \cite{Ecker-Huisken}, we can find $T>0$ and $C>0$ that make the lemma hold.

\end{proof}

Next, following the method of Lemma 6.5 in \cite{LN2021}, we prove a similar result for the case of the higher codimensions.

\begin{Lemma}\label{lemma_int_mean_curvature}
\begin{equation*}
    \underset{i\rightarrow\infty}{\lim}\,\dfrac{1}{\text{area}_h(D_i)}\int_{D_i}|H_h|^2dA_h=0,
\end{equation*}
where $H_h$ denotes the mean curvature of each disc $D_i$ in $(M,h)$.
\end{Lemma}

\begin{proof}
Suppose by contradiction that there exists $\epsilon>0$, such that after passing to a subsequence, and for $i\in\mathbb{N}$ large enough,
\begin{equation}\label{2epsilon}
    \dfrac{1}{\text{area}_h(D_i)}\int_{D_i}|H_h|^2dA_h> 2\epsilon.
\end{equation}
Under the mean curvature flow, the mean curvature satisfies the following evolution equation on the time interval $t\in[0,T]$ (see \cite{smoczyk2011mean}). \begin{align*}
    \nabla_{\frac{d}{dt}}|H_{h}(t)|^2= &\Delta|H_{h}(t)|^2-2|\nabla H_{h}(t)|^2+4\langle A_{h}^{jk}(t),H_{h}(t)\rangle_{h_t}\langle (A_{h})_{jk}(t),H_{h}(t)\rangle_{h_t}\\
    &+2(R_{h})_{jklm}(t)H_{h}^j(t)(F_i)_p^k(t) H_{h}^l(t)F_i^{mp}(t),
\end{align*}
where $H_{h}(t)$, $A_{h}(t)$ and $F_i(t)$ represent the mean curvature, second fundamental form and the immersion $F_i(t):D_i(t)\rightarrow M$, respectively. 

Using the result of Lemma \ref{short_time_existence}, we can pick a uniform constant $C_1>0$, such that for all sufficiently large $i\in\mathbb{N}$, and for any $t\in[0,T]$, 
\begin{equation*}
    {\frac{d}{dt}}\int_{D_i(t)} |H_h(t)|^2 dA_h\geq -C_1\text{area}_h(D_i(t))\geq -C_1\text{area}_h(D_i),
\end{equation*}
the latter inequality follows from the fact \begin{equation*}
    \frac{d}{dt}\text{area}_{h}(D_i(t))=\frac{1}{2}\int_{D_i(t)}tr\langle \frac{d}{dt}h_{jk},h^{jk}\rangle\,dA_h
    =-\int_{D_i(t)} |H_h(t)|^2dA_h
    \leq 0.
\end{equation*}

We can choose $T_1<\min\{\frac{\epsilon}{C_1},T\}$, by assumption (\ref{2epsilon}), for any $i\in\mathbb{N}$ and $t\in [0,T_1]$, \begin{equation*}
    \int_{D_i(t)} |H_h(t)|^2 dA_h\geq \epsilon\,\text{area}_h(D_i)\geq \epsilon\,\text{area}_h(D_i(t)).
\end{equation*}
Then we obtain \begin{equation*}
    \frac{d}{dt}\text{area}_{h}(D_i(t))
    =-\int_{D_i(t)} |H_h(t)|^2dA_h
    \leq -\epsilon\,\text{area}_h(D_i(t)).
\end{equation*}
Thus, for any sufficiently large $i\in\mathbb{N}$, \begin{equation*}
    \frac{\text{area}_h(\Pi_i)}{\text{area}_h(D_i)}\leq \frac{\text{area}_h(D_i(t))}{\text{area}_h(D_i)}\leq \frac{e^{-\epsilon T_1}\text{area}_h(D_i)}{\text{area}_h(D_i)}=e^{-\epsilon T_1}<1,
\end{equation*}
which  violates (\ref{limit_1}).
\end{proof}

Furthermore, arguing like Lemma 4.2 of \cite{Jiang}, we deduce the following result from Lemma \ref{lemma_int_mean_curvature}. 
For any round circle $c\subset\partial_\infty \mathbb{H}^n$, it has a dense $\pi_1(M)$-orbit in $\partial_\infty\mathbb{H}^n$. In addition, $c$ can be represented by
$\underset{i\rightarrow\infty}{\lim}\Lambda(\phi_i\Pi_i\phi_i^{-1})$, where $\phi_i\in \pi_1(M)$, and $\Lambda(\phi_i\Pi_i\phi_i^{-1})$ represents the limit set of $\phi_i\Pi_i\phi_i^{-1}$. Redefine $D_i$ by the lifts of $S_i$ to $\mathbb{H}^n$ preserved by $\phi_i\Pi_i\phi_i^{-1}$.
It has the property that \begin{equation*}
    \lim_{i\rightarrow\infty} \int_{D_i\cap B_{R_i}(0)}|H_h|^2 dA_h=0,\quad R_i\rightarrow\infty.
\end{equation*}
Note that after passing to a subsequence, $D_i$ converges to the totally geodesic disc $D(c)\subset\mathbb{H}^n$ that is asymptotic to $c$.
Therefore, the mean curvature $H_h$ vanishes on $D(c)$, namely, $D(c)$ is a minimal disc of $B^n$ with respect to the metric $h$. And since $c$ is chosen arbitrarily, every totally geodesic disc of $\mathbb{H}^n$ must be  minimal for $h$.

We apply the result below for surfaces in 3-manifolds, the proof can be found in \cite{LN2021}.
\begin{Lemma}
Every totally geodesic disc in $\mathbb{H}^3$ is minimal with respect to another metric $h$ if and only if for any geodesic $\gamma\subset\mathbb{H}^3$, the following function is a constant \begin{equation*}
    t\mapsto |h|_{h_0}^{-\frac{1}{2}}h\,(\gamma'(t),\gamma'(t)).
\end{equation*}
\end{Lemma}

Because of the ergodicity of the geodesic flow in $(M,h_0)$, we can choose a geodesic $\gamma$ of $M$ whose orbit is dense in the unit tangent bundle. Let $\Tilde{\gamma}$ be the lift of $\gamma$ to $\mathbb{H}^n$. $\Tilde{\gamma}$ must be contained in a hyperbolic 3-ball $B\approx\mathbb{H}^3$. Applying the previous lemma to the  geodesic $\Tilde{\gamma}$ and ambient manifold $B$, we conclude that \begin{equation*}
    t\mapsto |h|_B|_{h_0}^{-\frac{1}{2}}h|_B\,(\Tilde{\gamma}'(t),\Tilde{\gamma}'(t))
\end{equation*} 
is constant. So the projection $\gamma$ in $M$ also satisfies that 
\begin{equation*}
    t\mapsto |h|_{h_0}^{-\frac{1}{2}}h\,(\gamma'(t),\gamma'(t))
\end{equation*} 
is constant. Thus due to the density, there is a constant $c>0$, such that for any vector field $X$ of the unit tangent bundle of $M$, \begin{equation*}
    |h|_{h_0}^{-\frac{1}{2}}h(X,X)=c\,h_0(X,X)\quad\Longrightarrow \quad |h|_{h_0}^{-\frac{1}{2}}h=c\,h_0.
\end{equation*}
As a result, $h$ coincides with a multiple of $h_0$.

\bibliographystyle{plain} 
\bibliography{ref}   

\begin{thebibliography}{10}

\bibitem{Besse}
A.~L. Besse.
\newblock {\em Einstein Manifolds}, volume~10 of {\em Ergebnisse der Mathematik
  und ihrer Grenzgebiete}.
\newblock Springer Berlin, Heidelberg, 1987.

\bibitem{Besson-Courtois-Gallot}
G.~Besson, G.~Courtois, and S.~Gallot.
\newblock Volume et entropie minimale des espaces localement sym\'{e}triques.
\newblock {\em Inventiones mathematicae}, 103(2):417--446, 1991.

\bibitem{Besson-Courtois-Gallot-2}
G.~Besson, G.~Courtois, and S.~Gallot.
\newblock Entropies et rigidit\'{e}s des espaces localement sym\'{e}triques de
  courbure strictement n\'{e}gative.
\newblock {\em Geometric and Functional Analysis}, 5:731--799, 1995.

\bibitem{CMN}
D.~Calegari, F.~Marques, and A.~Neves.
\newblock Counting minimal surfaces in negatively curved 3-manifolds.
\newblock {\em Duke Mathematical Journal}, 2022.

\bibitem{Chau}
A.~Chau, J.~Chen, and W.~He.
\newblock Lagrangian mean curvature flow for entire lipschitz graphs.
\newblock {\em Calculus of Variations and Partial Differential Equations},
  44:199--220, 2012.

\bibitem{Ebin}
D.~G. Ebin.
\newblock On the space of {Riemannian} metrics.
\newblock {\em Bulletin of the American Mathematical Society},
  74(5):1001--1003, 1968.

\bibitem{Ecker-Huisken}
K.~Ecker and G.~Huisken.
\newblock Interior estimates for hypersurfaces moving by mean curvature.
\newblock {\em Inventiones mathematicae}, 105(3):547--570, 1991.

\bibitem{Gromov}
M.~Gromov.
\newblock Foliated plateau problem, part {I}: Minimal varieties.
\newblock {\em Geometric and Functional Analysis}, 1:14--79, 1991.

\bibitem{Hamenstadt}
U.~Hamenst{\"a}dt.
\newblock Incompressible surfaces in rank one locally symmetric spaces.
\newblock {\em Geometric and Functional Analysis}, 25, 02 2014.

\bibitem{Hamilton}
R.~S. Hamilton.
\newblock Heat equations in geometry.
\newblock {\em Lecture notes}, Honolulu, Hawaii, 1989.

\bibitem{Hasselblatt-Katok}
B.~Hasselblatt and A.~Katok.
\newblock {\em Chapter 1 Principal structures}, volume~1 of {\em Handbook of
  Dynamical Systems}, pages 1--203.
\newblock Elsevier Science, 2002.

\bibitem{Hempel}
J.~Hempel.
\newblock {\em 3-manifolds}.
\newblock Annals of Mathematics Studies 86, Princeton University Press, 1976.

\bibitem{Jiang}
R.~Jiang.
\newblock Counting essential minimal surfaces in closed negatively curved
  n-manifolds.
\newblock {\em arXiv:2108.01796}, 2021.

\bibitem{Kahn-Labourie-Mozes}
J.~Kahn, F.~Labourie, and S.~Mozes.
\newblock Surface groups in uniform lattices of some semi-simple groups.
\newblock {\em arXiv: Differential Geometry}, 2018.

\bibitem{Kahn-Markovic}
J.~Kahn and V.~Markovic.
\newblock Counting essential surfaces in a closed hyperbolic 3-manifold.
\newblock {\em Geom. Topol.}, 16:601--624, 2012.

\bibitem{Knieper}
G.~Knieper.
\newblock Volume growth, entropy and the geodesic stretch.
\newblock {\em Mathematical Research Letters}, 2:39--58, 1995.

\bibitem{Labourie}
F.~Labourie.
\newblock Asymptotic counting of minimal surfaces in hyperbolic manifolds
  [according to {C}alegari, {M}arques and {N}eves].
\newblock {\em arXiv:2203.09366}, 2022.

\bibitem{Lowe}
B.~Lowe.
\newblock Deformations of totally geodesic foliations and minimal surfaces in
  negatively curved 3-manifolds.
\newblock {\em Geom. Funct. Anal.}, 31:895--929, 2021.

\bibitem{LN2021}
B.~Lowe and A.~Neves.
\newblock Minimal surface entropy and average area ratio.
\newblock {\em arXiv:2110.09451}, 2021.

\bibitem{Ratner}
M.~Ratner.
\newblock Raghunathan's topological conjecture and distributions of unipotent
  flows.
\newblock {\em Duke Mathematical Journal}, 63(1):235--280, 1991.

\bibitem{SY}
R.~Schoen and S.~T. Yau.
\newblock Existence of incompressible minimal surfaces and the topology of
  three dimensional manifolds with non-negative scalar curvature.
\newblock {\em Annals of Mathematics}, 110(1):127--142, 1979.

\bibitem{Shah}
N.~A. Shah.
\newblock Closures of totally geodesic immersions in manifolds of constant
  negative curvature.
\newblock {\em Singapore: World Scientific}, 1991.

\bibitem{smoczyk2011mean}
K.~Smoczyk.
\newblock Mean curvature flow in higher codimension - introduction and survey.
\newblock In {\em Global Differential Geometry}, pages 231--274. Springer
  Berlin Heidelberg, 2012.

\bibitem{Uhlenbeck}
K.~Uhlenbeck.
\newblock {\em Closed minimal surfaces in hyperbolic 3-manifolds}, pages
  147--168.
\newblock Seminar On Minimal Submanifolds. (AM-103). Princeton University
  Press, 1983.

\bibitem{Wang1993}
S.~Wang.
\newblock The $\pi_1$-injectivity of self-maps of nonzero degree on
  3-manifolds.
\newblock {\em Mathematische Annalen}, 297(1):171--190, 1993.

\bibitem{Warner}
F.~W. Warner.
\newblock Extension of the {Rauch} comparison theorem to submanifolds.
\newblock {\em Transactions of the American Mathematical Society},
  122(2):341--356, 1966.

\bibitem{White}
B.~White.
\newblock A local regularity theorem for mean curvature flow.
\newblock {\em Annals of mathematics}, 161(3):1487--1519, 05 2005.

\end{thebibliography}
\end{document}